\documentclass[11pt,a4paper,twoside]{article}
\usepackage{amsmath}
\usepackage{mythm}
\usepackage{epsfig}
\def\ds{\displaystyle}

\newcommand{\Ftimes}{\mathbf{F}^\times}
\def\mymark#1{\markboth{Non-Archimedean Function Theory: #1}
{Non-Archimedean Function Theory: #1}}
\newtheorem{theorem}{Theorem}[subsection]
\newtheorem{lemma}[theorem]{Lemma}
\newtheorem{prop}[theorem]{Proposition}
\newtheorem{cor}[theorem]{Corollary}

\newtheorem{exercise}[theorem]{Exercise}
\newtheorem{statement}[theorem]{Statement}
\newtheorem*{definition}{Definition}

\theoremstyle{remark}
\newtheorem{remark}[theorem]{Remark}
\newtheorem*{remark*}{Remark}
\newtheorem*{example*}{Example}
\newtheorem{example}[theorem]{Example}

\begin{document}
\mymark{Analogs of Basic Complex Function Theory}
\pagestyle{myheadings}

\title{Lectures on Non-Archimedean Function Theory\\
Advanced School on $p$-Adic Analysis and Applications\\
The Abdus Salam International Centre for Theoretical Physics\\
Trieste, Italy}
\author{William Cherry\\
Department of Mathematics\\
University of North Texas}
\date{August 31 -- September 4, 2009}
\maketitle

{\Large
\begin{list}{}{\labelwidth+80pt \leftmargin+80pt
}
\item[\textbf{Lecture 1:}] Analogs of Basic Complex Function
Theory

\item[\textbf{Lecture 2:}] Valuation Polygons and 
a Poisson-Jensen Formula

\item[\textbf{Lecture 3:}] Non-Archimedean Value Distribution Theory

\item[\textbf{Lecture 4:}] Benedetto's Non-Archimedean Island Theorems

\end{list}

\vfil
}

\thispagestyle{empty}

\newpage
\thispagestyle{empty}
~
\newpage
This lecture series is an introduction to non-Archimedean function theory.
The audience is assumed to be familiar with non-Archimedean fields and
non-Archimedean absolute values, as well as to have had a
standard introductory course in complex function theory. 
A standard reference for
the later could be, for example, \cite{Ahlfors}.
No prior exposure
to non-Archimedean function theory is supposed.  Full details on the
basics of non-Archimedean absolute values and the
construction of $p$-adic number fields, the most important of the
non-Archimedean fields, can be found in \cite{Robert}.

\section{Analogs of Basic Complex Function Theory}

\subsection{Non-Archimedean Fields}
Let $A$ be a commutative ring.
A \textbf{non-Archimedean absolute value}
$|~|$ on $A$ is a function from $A$ to 
the non-negative real numbers $\mathbf{R}_{\ge0}$ satisfying the following
three properties:
\begin{quote}
\begin{description}
\item[\textbf{AV 1.}] $|a|=0$ if and only if $a=0;$
\item[\textbf{AV 2.}] $|ab|=|a|\cdot|b|$ for all $a,b\in A;$ and
\item[\textbf{AV 3.}] $|a+b| \le \max\{|a|,|b|\}$ for all $a,b\in A.$
\end{description}
\end{quote}

\begin{exercise}\label{sumeqmax} Prove that \textup{\textbf{AV 3}}
implies that if
$|a|\ne|b|,$ then $|a+b|=\max\{|a|,|b|\}.$
\end{exercise}

\begin{remark*} There are two geometric interpretations
of Exercise~\ref{sumeqmax}.
The first is that every triangle in a non-Archimedean world is isosceles.
The second is that every point inside a circle may serve as
a center of the circle. This also means that either two discs are
disjoint or one is contained inside the other.
\end{remark*}

\begin{exercise}\label{extendtofracfield} If $|~|$ is a non-Archimedean
absolute value on an integral domain $A,$ prove that $|~|$ extends uniquely
to the fraction field of $A.$
\end{exercise}

A pair $(\mathbf{F},|~|)$ consisting of a field $\mathbf{F}$ together with
a non-Archimedean absolute value $|~|$ on $\mathbf{F}$ will be referred to
as a \textbf{non-Archimedean field} and denoted simply by $\mathbf{F}$
for brevity. A sequence $a_n$ in a non-Archimedean field $\mathbf{F}$
is said to \textbf{converge} to an element $a$ in $\mathbf{F},$ if
for every $\varepsilon>0,$ there exists a natural number $N$ such that
for all natural numbers $n\ge N,$ we have \hbox{$|a_n-a|<\varepsilon.$}
A sequence $a_n$ in $\mathbf{F}$ is called a \textbf{Cauchy sequence}
if for every $\varepsilon>0,$ there exists a natural number $N$ such
that if $m$ and $n$ are both natural numbers $\ge N,$ then
$|a_n-a_m|<\varepsilon.$ As in elementary analysis, it is easy to see
that every convergent sequence is Cauchy. In general, not every Cauchy
sequence must converge.  However, if the non-Archimedean
field $\mathbf{F}$ is such that
every Cauchy sequence of elements in $\mathbf{F}$ converges, then
$\mathbf{F}$ is called \textbf{complete}.

\begin{exercise}\label{completion}
Let $\mathbf{F}$ be a non-Archimedean field.
Let $\overline{\mathbf{F}}$ be the set of Cauchy sequences in
$\mathbf{F}$ modulo the sequences which converge to $0.$ In other words,
define an equivalence relation on the set of Cauchy sequences in
$\mathbf{F}$ by defining two Cauchy sequences to be equivalent if their
difference is a sequence which converges to $0,$ and let
$\overline{\mathbf{F}}$ be the set of equivalence classes under this
equivalence relation.  Then, show $\overline{\mathbf{F}}$ is field,
that $|~|$ naturally extends to $\overline{\mathbf{F}},$ and that
$\overline{\mathbf{F}}$ is a complete non-Archimedean field, which we
call the \textbf{completion} of $\mathbf{F}.$
\end{exercise}

Given a field $\mathbf{F},$  we use $\mathbf{F}^\times$ to denote
$\mathbf{F}\setminus\{0\}.$
Given a non-Archimedean field $(\mathbf{F},|~|),$ the set
$|\Ftimes|=\{|x| : x \in\mathbf{F}^\times\}\subset\mathbf{R}_{>0}$
is a subgroup under multiplication of $\mathbf{R}_{>0}$ and is called
the \textbf{value group} of $\mathbf{F}.$ If $|\Ftimes|$ is discrete
in $\mathbf{R}_{>0},$ then $\mathbf{F}$ is called a 
\textbf{discretely valued} non-Archimedean field.

\smallskip
We now present some fundamental examples of non-Archimedean fields.

\subsubsection*{The Trivial Absolute Value}
Let $\mathbf{F}$ be any field. Define an absolute value $|~|,$ called
the \textbf{trivial absolute value}, 
on $\mathbf{F}$ by declaring that $|0|=0$ and $|x|=1$ for
all $x$ in $\mathbf{F}^\times.$ Clearly a sequence is Cauchy with
respect to the trivial absolute value if and only if it is eventually
constant, and hence convergent.  Thus any field can be made into
a complete non-Archimedean field by endowing it with the trivial absolute
value.

\subsubsection*{$p$-Adic Number Fields}
Consider the rational numbers $\mathbf{Q},$ and let $p$ be a prime number.
Then, any non-zero $x$ in $\mathbf{Q}$ can be written as
$$
	x=p^n\frac{a}{b},
$$
where $p$ does not divide $a$ or $b.$ If we define $|x|_p=p^{-n}$
and $|0|_p=0,$ then we easily see that $|~|_p$ is a non-Archimedean
absolute value on $\mathbf{Q}.$

\begin{exercise}\label{Qpseries}
Let $p$ be a prime number, let $n_0$ be an integer, and for each integer
$n\ge n_0,$ let $a_n$ be an integer between zero and $p-1,$ inclusive.
Show that sequence of partial sums
$$
	S_k = \sum_{n=n_0}^k a_n p^n
$$
is a Cauchy sequence in $(\mathbf{Q},|~|_p).$ Moreover, show that
$S_k$ converges in $\mathbf{Q}$ if and only if the $a_n$ are eventually
periodic, or in other words there exists integers $n_1$ and $t\ge1$
such that $a_{n+t}=a_n$ for all $n\ge n_1.$
\textit{Hint:} A solution can be found in 
\textup{\cite[\S I.5.3]{Robert}}.
\end{exercise}

We conclude from Exercise~\ref{Qpseries} that $\mathbf{Q}$ is not
complete with respect to $|~|_p,$ because, for example, the sequence
of partial sums
$$
	S_k = \sum_{n=0}^k p^{n^2}
$$
is Cauchy, but not convergent.  We denote by $\mathbf{Q}_p$
the completion of $\mathbf{Q}$ with respect to $|~|_p$ and call
this field the field of \textbf{$p$-adic numbers}. The closure
of the integers $\mathbf{Z}$ in $\mathbf{Q}_p$ is denoted by
$\mathbf{Z}_p,$  and elements of $\mathbf{Z}_p$ are called
\textbf{$p$-adic integers}.

\begin{exercise} Fix a prime number $p.$ Every non-zero
element $x$ in $\mathbf{Q}_p$
has a unique $p$-adic expansion of the form
$$
	x = \sum_{n=n_0}^\infty a_n p^n,
$$
where the $a_n$ are integers between $0$ and $p-1,$ $a_{n_0}\ne0,$
and $p^{-n_0}=|x|_p.$
\end{exercise}

\begin{exercise} Finite algebraic extensions of complete non-Archimedean
fields are again complete non-Archimedean fields.
Hint: See \textup{\cite[Ch.~XII]{Lang}.}
\end{exercise}

Finite algebraic extensions of $\mathbf{Q}_p$ are called
\textbf{$p$-adic number fields}.

\begin{exercise} No $p$-adic number field is algebraically closed.
Hint: Show that the value group of any finite extension of
$\mathbf{Q}_p$ must be discrete and hence cannot contain all the
$n$-th roots of $p$ for all $n.$
\end{exercise}

\begin{theorem}\label{CpThm}
The absolute value $|~|_p$ extends uniquely to the algebraic closure
$\mathbf{Q}_p^a$
of $\mathbf{Q}_p,$ which is not complete, but its completion
$\mathbf{C}_p$ remains algebraically closed.
\end{theorem}

The field $\mathbf{C}_p$ is called the \textbf{$p$-adic complex numbers}.
I will not discuss the proof of Theorem~\ref{CpThm} here.
See \cite[Ch.~III]{Robert} for a proof.

\subsubsection*{Positive Characteristic}

The following fields are important positive characteristic analogs
of the $p$-adic number fields.  Let $\mathbf{F}_q$ denote the finite
field of $q$ elements, where $q$ is a power of a prime. Let
$\mathbf{F}_q(T)$ denote the field of rational functions over 
$\mathbf{F}_q.$ In positive characteristic number theory, the
field $\mathbf{F}_q(T)$ plays the role of the rational numbers and
the polynomial ring $\mathbf{F}_q[T]$ plays the role of the integers.
The notation $|~|_\infty$ is often used to denote the unique
non-Archimedean absolute value on $\mathbf{F}_q(T)$ such
that $|T|_\infty=q.$ The completion of $\mathbf{F}_q(T)$ with respect
to $|~|_\infty$
is isomorphic to $\mathbf{F}_q((1/T)),$
the formal Laurent series ring in $1/T$ with coefficients in 
$\mathbf{F}_q.$ The complete non-Archimedean field
$(\mathbf{F}_q((1/T)),|~|_\infty)$ is a positive characteristic analog
of the $p$-adic number fields, namely the finite extensions of 
$\mathbf{Q}_p.$ As with the $p$-adic number fields, the absolute value
$|~|_\infty$ extends uniquely
to the algebraic closure of $\mathbf{F}_q((1/T)),$ and the completion
of $\mathbf{F}_q((1/T))^a$ remains algebraically closed, and is denoted
by $\mathbf{C}_\infty,$ or possibly
$\mathbf{C}_{p,\infty}$ if one wants to also emphasize the characteristic.
Hence $\mathbf{C}_\infty$ is a positive characteristic and non-Archimedean
analog of the complex numbers.

\medskip
These notes will discuss analysis over complete algebraically closed
non-Archimedean fields.  The most important examples of such fields
are the fields of $p$-adic complex numbers $\mathbf{C}_p$ and the
fields $\mathbf{C}_\infty$ introduced above. However, rarely is the
precise form of the field important, and henceforth $\mathbf{F}$
will denote simply a complete non-Archimedean field. Sometimes we may need to
assume that $\mathbf{F}$ has characteristic zero.

\subsection{Analytic and Meromorphic Functions}
Let $(\mathbf{F},|~|)$ be a complete, algebraically closed, non-Archimedean
field.

\begin{exercise}\label{easyconvergence}
A series $\sum a_n$ of elements of $\mathbf{F}$
converges if and only if $\lim\limits_{n\to\infty}|a_n|=0.$
\end{exercise}

Because of Exercise~\ref{easyconvergence}, there is no need in 
non-Archimedean analysis for any of the various convergence tests
one learns in freshmen calculus.

\smallskip
The formal power series ring $\mathbf{F}[[z]]$ in the variable $z$
with coefficients in $\mathbf{F}$ forms an integral domain with
addition and multiplication defined in the natural way.
Because of Exercise~\ref{easyconvergence}, an element
$$
	f(z)=\sum_{n=0}^\infty a_n z^n \in \mathbf{F}[[z]]
$$
is seen to converge at the point $z$ in $\mathbf{F}$ if
$$
	\lim_{n\to\infty} |a_n||z|^n=0.
$$
If a formal power series $f$ converges at $z,$ then clearly
$f$ converges at each $w$ with $|w|<|z|.$ Similarly, if
$f$ diverges at $z,$ then $f$ diverges at each $w$ with $|w|>|z|.$
We therefore define the \textbf{radius of convergence} $r_f$ of a
formal power series $f$ by
$$
	r_f = \sup\{|z| : f \textrm{~converges at~} z\}.
$$
One then has the usual Hadamard formula for the radius of convergence
of a formal power series we are familiar with from real or complex
analysis.

\begin{exercise}[Hadamard Formula]\label{radiusofconvergence}
$\ds r_f = \frac{1}{\ds\limsup_{n\to\infty}|a_n|^{1/n}}.$
\end{exercise}

It is also easy to see that radius of convergence behaves well
under addition and multiplication:

\begin{exercise}\label{radiusaddmult}
$\ds r_{f+g}\ge \min\{r_f,r_g\}$ and
$\ds r_{fg}\ge\min\{r_f,r_g\}.$
\end{exercise}

Define the \textbf{open} or \textbf{unbordered} \textbf{ball} of radius
$R$ by 
$$
	\mathbf{B}_{<R} = \{z \in \mathbf{F} : |z| < R\}.
$$
We also use the notation $\mathbf{B}_{<\infty}=\mathbf{F}$ to include
the case of all of $\mathbf{F}.$
The \textbf{closed} or \textbf{bordered} \textbf{ball} of radius $R<\infty$
is defined by
$$
	\mathbf{B}_{\le R} = \{z \in\mathbf{F} : |z|\le R\}.
$$
If $R>0,$ then both $\mathbf{B}_{<R}$ and $\mathbf{B}_{\le R}$
are both open and closed in the topology on $\mathbf{F}.$ Because of
this some people prefer the somewhat more cumbersome ``unbordered''
and ``bordered'' terminology. The ring of \textbf{analytic functions}
on $\mathbf{B}_{\le R},$ denoted $\mathcal{A}[R],$ is defined by
$$
	\mathcal{A}[R] = \left\{ \sum_{n=0}^\infty a_n z^n \in
	\mathbf{F}[[z]] : \lim_{n\to\infty} |a_n| R^n = 0\right\}.
$$
Similarly, the ring of analytic functions on $\mathbf{B}_{<R},$ denoted
$\mathcal{A}(R),$ is defined by
$$
	\mathcal{A}(R) = \left\{ \sum_{n=0}^\infty a_n z^n \in
	\mathbf{F}[[z]] : \lim_{n\to\infty} |a_n| r^n = 0
	\textrm{~for all~}r<R\right\}.
$$
Elements of $\mathcal{A}(\infty),$ \textit{i.e.} power series with
infinite radius of convergence, are called \textbf{entire} functions.

All this extends easily to convergent Laurent series.  Namely,
we can consider various types of bordered, unbordered, or semi-bordered
annuli:
\begin{align*}
	A[r_1,r_2] &= \{z \in \mathbf{F} : r_1 \le |z| \le r_2\} \quad
	&A(r_1,r_2] &=\{z \in \mathbf{F} : r_1 < |z| \le r_2\}\\
	A[r_1,r_2) &=\{z \in \mathbf{F} : r_1 \le |z| < r_2\} \quad
	&A(r_1,r_2) &= \{z \in \mathbf{F} : r_1 < |z| < r_2\},
\end{align*}
and the various rings of analytic functions on those spaces
\begin{align*}
	\mathcal{A}[r_1,r_2] &= \left\{\sum_{n=-\infty}^\infty a_n z^n :
	\lim_{|n|\to\infty} |a_n|r^n=0 \textrm{~for all~} r_1\le r \le r_2
	\right\}\\
	\mathcal{A}(r_1,r_2] &= \left\{\sum_{n=-\infty}^\infty a_n z^n :
	\lim_{|n|\to\infty} |a_n|r^n=0 \textrm{~for all~} r_1< r \le r_2
	\right\}\\
	\mathcal{A}[r_1,r_2) &= \left\{\sum_{n=-\infty}^\infty a_n z^n :
	\lim_{|n|\to\infty} |a_n|r^n=0 \textrm{~for all~} r_1\le r < r_2
	\right\}\\
	\mathcal{A}(r_1,r_2) &= \left\{\sum_{n=-\infty}^\infty a_n z^n :
	\lim_{|n|\to\infty} |a_n|r^n=0 \textrm{~for all~} r_1< r < r_2
	\right\}.
\end{align*}

Notice that all the above rings of analytic functions are integral
domains.  Elements of their fraction fields are called \textbf{meromorphic
functions}. Thus, I will use, for instance $\mathcal{M}(r_1,r_2]$
to denote the fraction field of $\mathcal{A}(r_1,r_2],$ which is the
field of meromorphic functions on $A(r_1,r_2].$

For the most part, I will leave a discussion of analytic 
and meromorphic functions on 
subsets of $\mathbf{F}$ more complicated than annuli to the other lecturers
in this school, and in particular I refer the reader to Berkovich's lectures.

\begin{remark}\label{trivfunc}If the absolute value $|~|$ on 
$\mathbf{F}$ is trivial, then $\mathcal{A}(1)$ is simply the
formal power series ring $\mathbf{F}[[z]]$ and $\mathcal{A}[1]$
is
the polynomial ring $\mathbf{F}[z].$
The ring of analytic functions on the annulus $A[1,1]$
are simply elements of $\mathbf{F}[z,z^{-1}].$
\end{remark}

\subsection{The Schnirelman Integral and an Analog of the Cauchy Integral
Formula}

If you think back to your first course on complex function theory,
probably nothing stands out as much as the Cauchy Integral Theorem
and the Cauchy Integral Formula. Thus, I feel it is most appropriate
to begin with a discussion of an integral introduced by
Schnirelman \cite{Shnirelman} that serves as an analog to the path
integral around a circle so commonplace in complex analysis and from
which analogs of many of the usual first theorems in complex analysis
can be derived. I point out, however, that the Schnirelman integral is
not used much anymore, and the consequences of this integral that
I will explain in this section can also be derived by the techniques
to be introduced in future lectures. My lectures here on the 
Schnirelman integral are based on \cite{Adams}, but I have changed the
definition to make a closer parallel with classical complex function
theory.

\subsubsection*{Definition and Basic Properties}

Consider the homomorphism from $\mathbf{Z}$ to $\mathbf{F}$
defined by sending an integer $n$ to $n\cdot 1$ in $\mathbf{F}.$
If $\mathbf{F}$ has characteristic zero, this homomorphism is 
injective, and otherwise its image is the prime field of $\mathbf{F}.$
When we write $|n|$ in this section,
by abuse of notation, we will mean the absolute
value of the image of $n$ in $\mathbf{F},$ even when $\mathbf{F}$
has positive characteristic, so for instance $|n|=0$ if $n$ is divisible
by the characteristic of $\mathbf{F}.$

\begin{exercise}\label{smallints}
The set of $n$ in $\mathbf{Z}$ such that $|n|<1$ forms a prime ideal
of $\mathbf{Z}.$
\end{exercise}

As a consequence of Exercise~\ref{smallints}, there are infinitely many
positive integers $n$ such that $|n|=1.$

\begin{definition} For an integer $n\ge1$ such that $|n|=1,$
denote by $\xi_1^{(n)},\dots,\xi_n^{(n)}$ the $n$ $n$-th roots of
unity in $\mathbf{F}.$
Given $a$ and $r$ in $\mathbf{F}$
and given a function $f$ such that $f$ is defined at all points
of the form $a+r\xi_k^{(n)}$ for all $n\ge1$ with $|n|=1$ and all
$1\le k \le n,$ define
$$
	\int\limits_{|z-a|=|r|}\!\!\!\!\!\!f(z)dz = \lim_{\begin{array}{cc}
\noalign{\vskip -2pt}
	\scriptstyle n\to\infty\\
\noalign{\vskip -4pt}
	\scriptstyle |n|=1\end{array}} \frac{r}{n}\sum_{k=1}^n
	f\left(a+r\xi_k^{(n)}\right)\xi_k^{(n)},
$$
provided the limit on the right exists.
\end{definition}
The integral in the above definition is called the
\textbf{Schnirelman integral,} and if it exists, the function $f$
is called \textbf{Schnirelman integrable} on the discrete
circle $|z-a|=|r|.$ It is clear from the definition that the
Schnirelman integral satisfies the usual linearity properties we
expect from an integral.

\smallskip
\textbf{Caution!} The above definition for the Schnirelman integral
is non-standard.  Also, it could be that
$$
	\int\limits_{|z-a|=|r_1|}\!\!\!\!\!\!f(z)dz \;\;\;\ne
	\int\limits_{|z-a|=|r_2|}\!\!\!\!\!\!f(z)dz,
$$
or even that one of the above integrals exists and the other does not,
even if $|r_1|=|r_2|.$ This is not well-reflected in my choice of notation.

\begin{prop}\label{SchnirelmanEstimate}
If $\!\!\!\!\!\!\ds\int\limits_{|z-a|=|r|}\!\!\!\!\!\!f(z)dz$ exists, then
$$
	\left|\,\,\int\limits_{|z-a|=|r|}\!\!\!\!\!\!f(z)dz\right|
	\le |r|\max_{|z-a|=|r|} |f(z)|,
$$
provided the right hand side is well-defined.
\end{prop}

\begin{proof} Trivial, noting that $|\xi_k^{(n)}|=1.$
\end{proof}

\begin{prop}\label{SchnirelmanUnifConv}
If $\sum f_j$ converges uniformly on $|z-a|=|r|$ to $f$
and if each $f_j$ is Schnirelman integrable on $|z-a|=|r|,$
then $f$ is Schnirelman integrable on $|z-a|=|r|$ and
$$
	\int\limits_{|z-a|=|r|}\!\!\!\!\!\!f(z)dz = \sum_j
	\int\limits_{|z-a|=|r|}\!\!\!\!\!\!f_j(z)dz.
$$
\end{prop}

\begin{proof} Let $\varepsilon > 0.$ By the hypothesis of uniform
convergence of the sum, we have
$$
	\left|f(z)-\sum_{j=0}^Jf_j(z)\right|<\frac{\varepsilon}{3|r|}
$$
for all sufficiently large $J$ and for all $z$ such that $|z-a|=|r|.$
Hence, for any
$n$ such that $|n|=1,$
$$
	\left|\frac{r}{n}\sum_{k=1}^nf\left(a+r\xi_k^{(n)}\right)\xi_k^{(n)}
	\;-\;\frac{r}{n}\sum_{k=1}^n\sum_{j=0}^J
	f_j\left(a+r\xi_k^{(n)}\right)\xi_k^{(n)}\right|<\frac{\varepsilon}{3},
$$
for all sufficiently large $J.$
Since $|f_j|$ tends uniformly to zero on $|z-a|=|r|,$
we have by Proposition~\ref{SchnirelmanEstimate}, that
$$
	\left|
	\sum_{j=0}^J \;\;\int\limits_{|z-a|=|r|}\!\!\!\!\!\!f_j(z)dz
	\;-\;
	\sum_{j=0}^\infty\;\;\int\limits_{|z-a|=|r|}\!\!\!\!\!\!f_j(z)dz
	\right|<\frac{\varepsilon}{3},
$$
also for all sufficiently large $J.$ Fix $J$ sufficiently large
that the above inequalities hold.
By the integrability of the $f_j,$
there exists an $N$ such that if $n\ge N$ and $|n|=1,$ then
$$
	\left | \frac{r}{n}\sum_{k=1}^n\sum_{j=0}^J
	f_j\left(a+r\xi_k^{(n)}\right)\xi_k^{(n)}
	\;-\;
	\sum_{j=0}^J \;\;\int\limits_{|z-a|=|r|}\!\!\!\!\!\!f_j(z)dz\right|
	< \frac{\varepsilon}{3}.
$$
Hence, if $n\ge N$ and $|n|=1,$ then we have
$$
	\left|\frac{r}{n}\sum_{k=1}^nf\left(a+r\xi_k^{(n)}\right)
	\xi_k^{(n)}\;-\;
	\sum_{j=0}^\infty
		\;\;\int\limits_{|z-a|=|r|}\!\!\!\!\!\!f_j(z)dz\right|
	< \varepsilon.\qedhere
$$
\end{proof}

\subsubsection*{Cauchy Integral Theorem and Cauchy Integral Formula}

\begin{lemma}\label{nthroot} Let $1\le |j| < n$ be integers
\textup{[}here $|j|$
denotes the usual Archimedean absolute value of the index $j$\textup{]}. Then
$$
	\sum_{k=1}^n \left(\xi_k^{(n)}\right)^j=0.
$$
\end{lemma}

\begin{proof}
Since $\{\xi_1^{(n)},\dots,\xi_n^{(n)}\}=\{(\xi_1^{(n)})^{-1},
\dots,(\xi_n^{(n)})^{-1}\},$ it suffices to consider $j$ positive.
Let $x_1,\dots,x_n$ be variables.  Then,
$x_1^j+\dots+x_n^j$ is a polynomial in the elementary symmetric
functions $\sigma_1(x_1,\dots,x_n),\dots,\sigma_j(x_1,\dots,x_n)$
with no constant term.  Since $\sigma_i(\xi_1^{(n)},\dots,\xi_n^{(n)})=0$
for $1\le i < n,$ the lemma follows.
\end{proof}

\begin{theorem}[Cauchy Integral Theorem]\label{CauchyThm}
Let $\mathbf{B}_{\le R}(a)=\{z \in \mathbf{F} : |z-a|\le R\}$
denote the closed ball of radius $R$ centered at $a.$ Let $f$
be analytic on $\mathbf{B}_{\le R}(a).$ Let $r\in\mathbf{F}$
with $|r|=R.$ Then, $f$ is Schnirelman integrable on $|z-a|=|r|$ and
$$
	\int\limits_{|z-a|=|r|}\!\!\!\!\!\!f(z)dz=0.
$$
\end{theorem}

\begin{proof}
Without loss of generality, we may assume $a=0.$ By linearity
and Proposition~\ref{SchnirelmanUnifConv}, it
suffices to show the theorem for $f(z)=z^j$ for $j\ge0.$
The theorem then follows from Lemma~\ref{nthroot} since the 
expression inside the limit defining the Schnirelman integral
vanishes as soon as \hbox{$n\ge j+2.$}
\end{proof}

Given a formal power series $f(z)=\sum a_j z^j,$ define the
$k$-th Hasse derivative of $f$ by 
$$
	D^kf(z)=\sum_{j=k}^\infty a_j\binom{j}{k}z^{j-k}.
$$
Observe that in characteristic zero, the Hasse derivative
$D^kf$ is simply $f^{(k)}/k!.$ Hasse derivatives are more useful
than ordinary derivatives in positive characteristic and have
similar properties.

\begin{theorem}[Cauchy Integral Formula]\label{CauchyFmla}
Let $f$ be analytic on $\mathbf{B}_{\le R}(a),$
let $r\in\mathbf{F}$
with $|r|=R,$ let $w\in\mathbf{F}$ with $|w-a|\ne R,$ and let $n\ge0.$
Then
$$
	\int\limits_{|z-a|=|r|}\frac{f(z)}{(z-w)^{(n+1)}}\,dz
	=\left\{\begin{array}{ll}
	D^nf(w) &\textrm{if~}|w-a|<R\\
	0&\textrm{if~}|w-a|>R.
	\end{array}\right.
$$
\end{theorem}

\begin{proof}
From the definition and from Lemma~\ref{nthroot},
we see that if $k$ is an integer, then
\begin{equation}\label{CauchyEqn}
	\int\limits_{|z|=|r|}z^kdz=\left\{\begin{array}{ll}
	1&\textrm{if~}k=-1\\
	0&\textrm{otherwise}.
	\end{array}\right.
\end{equation}
Without loss of generality, we consider $a=0,$ and write $f(z)$
as a power series,
$$
	f(z)=\sum_{j=0}^\infty a_j z^j.
$$
If $w=0,$ then the theorem follows from~(\ref{CauchyEqn})
and Proposition~\ref{SchnirelmanUnifConv}. If $0<|w|<R,$ then
$$
	\sum_{k=n}^\infty\binom{k}{n}\left(\frac{w}{z}\right)^{k-n}
$$
converges uniformly to $1/(1-w/z)^{n+1}$ on $|z|=R.$ Hence,
$$
	\int\limits_{|z|=|r|}\frac{f(z)}{(z-w)^{n+1}}dz
	= \sum_{k=n}^\infty\sum_{j=0}^\infty
	w^{k-n}\int\limits_{|z|=|r|}a_j\binom{k}{n}z^{j-k-1}dz
	= \sum_{k=n}^\infty a_k w^{k-n}\binom{k}{n}=D^nf(w),
$$
where the second to last equality follows from~(\ref{CauchyEqn}).
If $|w|>R,$ we use a similar argument with 
$$
	\frac{1}{(z-w)^{n+1}}=\frac{(-1)^{n+1}}{w^{n+1}}\sum_{k=n}^\infty
	\binom{k}{n}\left(\frac{z}{w}\right)^k
$$
to conclude that the integral is zero, since there will be no negative
powers of $z.$
\end{proof}

\begin{theorem}[Residue Theorem]\label{ResidueThm}
Let $a$ and $r$ be elements of $\mathbf{F}.$ 
Let $f(z)$ be analytic on $|z-a|\le |r|,$ let $P(z)$ be a polynomial
with no zeros on $|z-a|=|r|,$ and let $R(z)=f(z)/P(z).$ Then,
$$
	\int\limits_{|z-a|=|r|}\!\!\!\!\!\!R(z)dz=
	\sum_{|b-a|<r}\mathrm{Res}(R,b).
$$
\end{theorem}

\begin{proof}Take the partial fraction expansion of $R,$
$$
	R(z)=g(z)+\frac{A_{1,1}}{z-b_1}+\dots+\frac{A_{1,m_1}}{(z-b_1)^{m_1}}
	+\dots+\frac{A_{n,1}}{z-b_n}+\dots+\frac{A_{n,m_n}}{(z-b_n)^{m_n}},
$$
where $g$ is analytic on $|z-a|\le|r|$
and apply Theorem~\ref{CauchyFmla}.
\end{proof}

\subsection{Consequences of the Cauchy Integral Formula}

\subsubsection*{Maximum Modulus Principle}

\begin{theorem}[Maximum Modulus Principle]\label{maxmod}
Let $r$ and $a$ be in $\mathbf{F},$ and let $f$ be analytic
on $|z-a|\le |r|.$ Then,
$$
	\max_{|z-a|\le|r|}|f(z)|=\max_{|z-a|=|r|}|f(z)|.
$$
\end{theorem}

\begin{proof}
Let $w$ be in $\mathbf{F}$ such that $|w-a|<|r|.$ Then, by
Theorem~\ref{CauchyFmla} and Proposition~\ref{SchnirelmanEstimate},
$$
	|f(w)| = \left|\,\int\limits_{|z-a|=|r|}\frac{f(z)}{z-w}dz\right|
	\le |r|\max_{|z-a|=|r|}\frac{|f(z)|}{|z-w|}.
$$
Now, $|z-w|=|(z-a)-(w-a)|=|z-a|=|r|$ by Exercise~\ref{sumeqmax}, and hence
the theorem follows.
\end{proof}

\begin{remark*}
In complex analysis, there is a stronger form of the maximum modulus
principle.  Namely if $f$ is analytic on $|z-a|=r$ and if $f$ attains its
maximum in the interior, then $f$ must be constant.  This is easily
seen to be false in non-Archimedean function theory. Indeed, consider
$|c|>1$ and let $f(z)=z+c.$ Then, $|f(z)|=|c|$ for all $|z|\le1,$
yet $f$ is not constant.  Notice, however, that $|f|$ is constant on
$|z|\le1.$
\end{remark*}

A variation on the maximum principle is

\begin{prop}\label{DerivEstimate}
Let $r$ and $a$ be in $\mathbf{F},$ and let $f$ be analytic
on $|z-a|\le |r|.$ Then,
$$
	|D^nf(w)|\le\frac{\ds\max_{|z-a|=|r|}|f(z)|}{|r|^n}
$$
for all $w$ in $\mathbf{F}$ such that $|w-a|<|r|$ and all integers
$n\ge0.$
\end{prop}

\begin{proof}
Fix $w$ in $\mathbf{F}$ such that $|w-a|<|r|.$ Then, by
Theorem~\ref{CauchyFmla} and Proposition~\ref{SchnirelmanEstimate},
as in the proof of the maximum principle, we have
$$
	|D^nf(w)| = \left|\;\int\limits_{|z-a|=|r|}\frac{f(z)}{(z-w)^{n+1}}
	dz\;\right|
	\le |r|\max_{|z-a|=|r|}\frac{|f(z)|}{|z-w|^{n+1}}.
$$
Again, $|z-w|=|(z-a)-(w-a)|=|z-a|=|r|$ by Exercise~\ref{sumeqmax}, and hence
the proposition follows.
\end{proof}

\begin{prop}\label{ChangeCenter} Let $a$ in $\mathbf{F}$
and let $f(z)=\sum_{j=0}^\infty c_j (z-a)^j$ be an analytic function with
radius of convergence $0<R\le\infty.$ Let $b$ be an element of
$\mathbf{F}$ such that $|b-a|<R.$ Then,
$$
	f(z)=\sum_{n=0}^\infty D^nf(b)(z-b)^n,
$$
and this second series has radius of convergence $R$ as well.
\end{prop}

\begin{proof}
The case of trivial absolute value, in which case both series are
polynomials or formal power series, is clear.
Hence, we assume the absolute value on $\mathbf{F}$
is non-trivial. Fix $z$ in $\mathbf{F}$ such that $|z-b| < R.$
Then $|z-a|=|(z-b)-(b-a)|<R,$ and we can find $r$ in $\mathbf{F}$ such
that
$$
	\max\{|z-b|,|b-a|\}<|r|<R.
$$
By Proposition~\ref{DerivEstimate},
$$
	|D^nf(b)(z-b)^n|\le M\left|\frac{z-b}{r}\right|^n,
\qquad\textrm{where~}M=\max_{|z-a|=|r|}|f(z)|.
$$
As $|(z-b)/r|<1,$ we see the series 
$$
	\sum_{n=0}^\infty D^nf(b)(z-b)^n
$$
converges at $z.$
Hence the radius of convergence of this series is at least $R.$
By symmetry, it is at most $R.$

\smallskip
Once we have convergence, we easily see
\begin{align*}
	\sum_{n=0}^\infty D^nf(b)(z-b)^n
	&=\sum_{n=0}^\infty\sum_{j=n}^\infty c_j
		\binom{j}{n}(b-a)^{j-n}(z-b)^n\\
	&=\sum_{j=0}^\infty c_j
		\sum_{n=0}^j\binom{j}{n}(b-a)^{j-n}(z-b)^n\\
	&=\sum_{j=0}^\infty c_j(z-b+b-a)^j=
	\sum_{j=0}^\infty c_j(z-a)^j=f(z).\qedhere
\end{align*}
\end{proof}

\subsubsection*{Identity Theorem}

\begin{theorem}[Identity Theorem]\label{IdentityTheorem}
Let $f$ be analytic on $\mathbf{B}_{\le R}(a)$ and let
$z_1,z_2,z_3,\dots$ be points in $\mathbf{B}_{\le R}(a)$ such that
$z_0$ is an accumulation point of the $z_k$ in $\mathbf{B}_{\le R}(a).$
If $f(z_k)=0$ for all $k\ge1,$ then $f\equiv 0.$
\end{theorem}

\begin{proof}
By Proposition~\ref{ChangeCenter}, we can expand $f$ as a power series
about $z_0,$
$$
	f(z)=\sum_{j=0}^\infty a_j(z-z_0)^j,
$$
and this series will also have radius of convergence $R.$
If not all the $a_j$ are zero, let $j_0$ be the smallest index
such that $a_{j_0}\ne0.$ Then, $f(z)/(z-z_0)^{j_0}$ is continuous and
non-zero for $|z-z_0|$ small. This contradicts the hypotheses
that $f(z_k)=0$ and $z_k$ accumulates at $z_0.$
\end{proof}

\subsubsection*{Liouville's Theorem}

\begin{theorem}[Liouville's Theorem]\label{LiouvilleThm}
A bounded entire function must be constant. Moreover, if
$$
	\frac{|f(z)|}{|z|^d}
$$
remains bounded as $|z|\to\infty,$ then $f$ must be a polynomial
of degree at most $d.$
\end{theorem}

\begin{proof}
Write $f$ as a power series
$$
	\sum_{j=0}^\infty a_j z^j.
$$
If $f$ is not a polynomial of degree at most $d,$ then there is
some coefficient $a_{j_0}\ne0$ with $j_0>d.$ By Theorem~\ref{CauchyFmla}
$$
	a_{j_0}=\int\limits_{|z|=|r|}\frac{f(z)}{z^{j_0+1}}dz.
$$
Then, by Proposition~\ref{SchnirelmanEstimate},
$$
	|a_{j_0}|\le\frac{\ds\max_{|z|=|r|}|f(z)|}{|r|^{j_0}}.
$$
Because $j_0>d,$ the right-hand side tends to zero as $|r|\to\infty,$
contradicting $a_{j_0}\ne0.$
\end{proof}

\subsubsection*{The Schwarz Lemma}

\begin{theorem}[Schwarz Lemma]
Let $f$ be analytic on $\mathbf{B}_{<1}$ such that the image of $f$
is contained in $\mathbf{B}_{\le 1}$ and such that $f(0)=0.$
Then, $|f(z)|\le|z|$ for all $|z|<1$ and $|f'(0)|\le1.$
\end{theorem}

\begin{proof}
If the absolute value on $\mathbf{F}$ is trivial, then so is the
statement of the theorem.  We therefore assume that the absolute value
is non-trivial. Hence, choose a sequence $r_n$ with $|r_n|<1$
and such that $|r_n|\to1.$ By the assumption that $f(0)=0,$
the function $g(z)=f(z)/z$ is also analytic on $\mathbf{B}_{<1}$
and by the maximum modulus principle, for $|z|<|r_n|,$ we have
$$
	|g(z)| \le \frac{1}{|r_n|}.
$$
Taking the limit as $|r_n|\to1$ completes the proof of the theorem.
\end{proof}

\begin{remark*}In the complex Schwarz Lemma, one has the additional statement
that if equality holds at some point, then $f(z)=\alpha z$ for some
$|\alpha|=1.$ This is easily seen to be false in the non-Archimedean
case by considering, for example, $f(z)=z(1+z).$
\end{remark*}

\begin{cor}[Schwarz-Pick]\label{SchwarzPick}
Let $f$ be analytic on $\mathbf{B}_{<1}$ such that the image of
$f$ is contained in $\mathbf{B}_{\le 1}.$ Then, for all $z$ and $w$
in $\mathbf{B}_{<1},$ we have $|f(z)-f(w)|\le|z-w|.$
\end{cor}

\begin{proof}
Fix $z$ and $w$ in $\mathbf{B}_{<1}.$
Consider $g(\zeta)=f(\zeta+w)-f(w).$ Then, if $|\zeta|<1,$
we have $|\zeta+w|<1,$ and so
$$
	|g(\zeta)|=|f(\zeta+w)-f(w)| \le 1,
$$
and hence $g$ satisfies the hypotheses of the theorem.
Now choosing $\zeta=z-w$ gives the corollary.
\end{proof}

\begin{remark*} The complex Schwarz-Pick lemma says that holomorphic
self-maps of the unit disc are distance non-increasing in the hyperbolic
metric.  Corollary~\ref{SchwarzPick} says that analytic self-maps
of a non-Archimedean disc are distance non-increasing in the standard
non-Archimedean metric.
\end{remark*}

\subsection{Morera's Theorem?}

I will conclude my first lecture by pointing out that
there is no converse of the Cauchy Integral Theorem for Schnirelman
integrals, that is no analog of Morera's Theorem.

\begin{example*}
Let $\mathbf{F}$ be a complete non-Archimedean field
which contains transcendental numbers and
consider the function $f$ which is 1 at all algebraic numbers in $\mathbf{F}$
and $0$ at all transcendental numbers of $\mathbf{F}.$ Clearly $f$ is
not given by a power series. On the other hand, let $a$ and $r\ne0$ be
elements of $f.$ Fix a positive integer $n$ with $|n|=1.$
Suppose that two of the numbers $a+r\xi_i^{(n)}$ and
$a+r\xi_j^{(n)}$ are algebraic. Then, $r(\xi_i^{(n)}-\xi_j^{(n)})$ is
algebraic and non-zero;  hence, $r$ is algebraic.
Therefore $a$ is also algebraic.
Thus, given $r$ and $a,$ one of three things can happen: $a+r\xi_i^{(n)}$
is algebraic for all $i=1,\dots,n;$
$a+r\xi_i^{(n)}$
is transcendental for all $i=1,\dots,n;$ or
$a+r\xi_i^{(n)}$
is algebraic for exactly one $i=1,\dots,n.$ In any of these cases,
we see that the sum defining the Schnirelman integral of $f$ tends to zero as
$n\to\infty.$
\end{example*}

\subsection{Concluding Remarks}

The Schnirelman integral seems never to have been a widely known
technique and is not often used in non-Archimedean function theory.
The methods I will discuss in the next lecture can be used to prove
the results of this section, and in fact some stronger results.
However,  the
Schnirelman integral is a
useful construct to have in one's bag of tricks because, in particular,
it often allows one to convert standard complex analytic proofs to
the non-Archimedean setting.  For example, in the 1960's, Adams \cite{Adams}
made extensive use of the Schnirelman integral to prove $p$-adic
versions of the Gelfond-Schneider-Lang transcendence machinery, and 
many of his proofs follow the same general outline of their complex
counterparts.

\section{Valuation Polygons and a Poisson-Jensen Formula}
\mymark{Valuation Polygons and a Poisson-Jensen Formula}

In the previous lecture, Schnirelman integrals were introduced so that
non-Archimedean analogs of familiar results from classical complex
function theory could be proved in a manner reminiscent of the proofs
most familiar from an introductory course in complex analysis.
Although the viewpoint of the previous lecture is nice in that it emphasizes
similarity between complex and non-Archimedean function theory, in fact,
there are many differences between non-Archimedean function theory and
its classical complex counterpart. What has been, in practice, a more
useful tool than the Schnirelman integral is a set of techniques for
determining the locations of non-Archimedean zeros of power series
that goes by either the name of the the valuation polygon or the 
Newton polygon. This is a powerful technique available in non-Archimedean
analysis that has no exact parallel in complex analysis.
Mastering these techniques is essential for developing non-Archimedean
analysis, and these techniques can generally be used in place of the
Cauchy Integral Formula and often give stronger results.
Takashi Harase, in his review \cite{Harase} of one of Dwork's
last papers \cite{Dwork}, writes:
\begin{quote}
``The author uses the Newton polygon at the places where classical
analysts use Cauchy's integral formula. He could be called the magician
of the Newton polygon.''
\end{quote}

\medskip
Much of the text of this lecture is taken from \cite{CherryWang}.

\subsection{The Residue Class Field}

Non-Archimedean fields have some structure not available in the complex
numbers or other Archimedean fields. Observe that the set
$$
	\mathcal{O}=\{z \in \mathbf{F} : |z| \le 1\}
$$
forms a subring of $\mathbf{F}.$ The subring $\mathcal{O}$ is called
the \textbf{ring of integers} of $\mathbf{F}.$ Denote by
$$
	M=\{z \in \mathbf{F} : |z| < 1 \}.
$$

\begin{exercise}Show that $M$ is the unique maximal ideal of $\mathcal{O}.$
\end{exercise}

\newcommand{\Ftilde}{\widetilde{\mathbf{F}}}
\newcommand{\F}{\mathbf{F}}
\renewcommand{\O}{\mathcal{O}}
\newcommand{\R}{\mathbf{R}}
\newcommand{\Z}{\mathbf{Z}}
A ring with a unique maximal idea is called a \textbf{local ring}.
We denote by $\Ftilde$ the field $\mathcal{O}/M,$
which is called the \textbf{residue class field of} $\mathbf{F}.$
Given an element $a$ in $\mathcal{O},$ we denote by $\tilde a$
its image in $\Ftilde.$

\begin{exercise}If $\mathbf{F}$ is algebraically closed, then
$\widetilde{\mathbf{F}}$ is too.
\end{exercise}

\subsection{Non-Archimedean Absolute Values on Rings of Analytic Functions}

Let $r_1\le r_2$ and consider the ring $\mathcal{A}[r_1,r_2]$
of analytic functions on the bordered annulus $A[r_1,r_2].$
The elements of $\mathcal{A}[r_1,r_2]$ are Laurent series of the
form
$$
	f(z)=\sum_{n=-\infty}^\infty c_n z^n
	\textrm{~such that~}\lim_{|n|\to\infty} |c_n|r^n=0
	\textrm{~for all~}r_1\le r \le r_2.
$$

\begin{remark*}
Unlike in the previous lecture, here $|n|$ for an integer index $n$
refers to the usual Archimedean absolute value of $n,$ and not
the absolute value of the image of $n$ in $\mathbf{F}.$
Also, in the previous lecture we tended to use $r$ to denote
an element of $\mathbf{F}.$ In this section, $r$ will denote a
non-negative real number.
\end{remark*}

For each $r$ between $r_1$ and $r_2,$ we define
$$
        \fbox{$|f|_r = \sup |c_n| r^n.$}
$$

\begin{remark*}Observe that for fixed $f,$ we easily see that
$|~|_r$ is a continuous function of $r.$ If $r_1=0,$ so that no negative
powers appear in the series expansion for $f,$ we also see that
$|f|_r$ is non-decreasing in $r.$
\end{remark*}

Next, we highlight one technical point.
\begin{remark*}  Recall that
we use $|\Ftimes|$ to denote the value group
$$
	|\Ftimes| = \{|a| \in \R_{>0} : a \in \F\}.
$$
It could be that $|\Ftimes|\ne\R_{>0},$ and this often creates some 
technicalities in non-Archimedean proofs where one has to consider
separately the cases $r\in|\Ftimes|$ and $r\not\in|\Ftimes|.$
\end{remark*}

We will see in Proposition~\ref{absrprop} below that 
$|~|_r$ is in fact a non-Archimedean absolute value if $r>0,$
but first we state another version
of the non-Archimedean maximum modulus principle.

\begin{prop}[Maximum Modulus Principle]
\label{MaxModProp}
If $f$ is analytic on $A[r_1,r_2]$ and $z_0$ is in $ A[r_1,r_2],$ then
$$
        |f(z_0)| \le |f|_{|z_0|}.
$$
Moreover, if $|z|=|z_0|$ and $|f(z)|<|f|_{|z_0|},$ then
$\widetilde{\frac{z}{z_0}}$ is one of at most finitely many residue
classes in $\Ftilde.$
\end{prop}

\begin{proof}
That $|f(z_0)| \le |f|_{|z_0|}$ follows immediately from the 
non-Archimedean triangle inequality, so we need to show that
equality holds outside at most finitely many residue classes.
Write 
$$
	f(z)=\sum_{n=-\infty}^\infty c_n z^n,
$$
and let
$c$ be an element of $\F$ such that
\hbox{$|c|=|f|_{|z_0|}.$} Let
$$
        g(z) = \sum_{n=-\infty}^\infty b_n z^n, \qquad\textrm{where~}
        b_n = \frac{c_n z_0^n}{c}.
$$
Note that $\sup |b_n|=1,$ and in particular $g$ has coefficients in 
$\O.$  Let 
$$
        \tilde g(z) = \sum_{n=-\infty}^\infty \tilde b_n z^n,
$$
and note that $\tilde g$ is not identically zero.  Then, 
directly from the definitions, if $|z|=|z_0|$ and
\hbox{$|f(z)|<|f|_{|z_0|},$} then 
$$
        \tilde g \left(\widetilde{\frac{z}{z_0}}\right)=0.
$$
Because 
$$
        \lim_{|n|\to\infty}|b_n|=0,
$$
$\tilde g$ has only finitely many non-zero coefficients, and
hence, there are only finitely many possibilities
for $\widetilde{\frac{z}{z_0}}.$
\end{proof}

\begin{cor}
\label{BoundedCor}
If $f$ is analytic on $A[r_1,r_2]$ and the set
\hbox{$\{|f|_r : r_1\le r \le r_2\}$} is bounded in $\R,$
then $f$ is bounded on $A[r_1,r_2].$
\end{cor}

\begin{prop}
\label{absrprop}
If $f$ and $g$ are analytic functions on $A[r_1,r_2],$ then
$$\begin{array}{rcll}
        |f+g|_r &\le& \max \{|f|_r, |g|_r\} \qquad&\textnormal{[non-Archimedean
triangle inequality]}\\
\noalign{\vspace{3pt}}
        |fg|_r &=& |f|_r |g|_r \qquad&\textnormal{[multiplicativity]}
\end{array}$$
for $r_1\le r \le r_2.$
\end{prop}

\begin{proof}
The triangle inequality is clear.  If $r\in|\Ftimes|,$
then multiplicativity follows from
Proposition~\ref{MaxModProp}. For $r\not\in|\Ftimes|,$ the case when
$|~|$ is trivial is left as an exercise for the reader. If $|~|$ is
non-trivial, we can find a sequence $r_n\in|\Ftimes|$ such that $r_n\to r,$
and then the proposition follows from the continuity of $|~|_r$ in $r.$
\end{proof}

Proposition~\ref{absrprop} says that 
each $|~|_r$ is a non-Archimedean absolute value on the ring of functions
analytic on an annulus containing $\{z : |z|=r\},$  provided $r>0.$

\begin{prop}
\label{CompleteProp}
The ring of analytic functions on $A[r_1,r_2]$ is complete
with respect to the norm 
$$
	|f|_{\sup} = \sup_{r_1\le r \le r_2}|f|_r.
$$
\end{prop}

\begin{proof}
Let
$$
	f_n(z) = \sum_m a_{n,m}z^m
$$
be a Cauchy sequence.  Then, for $n$ and $n'$ sufficiently large,
we have
\begin{equation}\label{cauchysupeqn}
	\varepsilon > |f_n-f_{n'}|_{\sup} = \sup_r \sup_m 
	|a_{n,m}-a_{n',m}|r^m.
\end{equation}
This implies that for each fixed $m,$ the sequence of coefficients
$a_{n,m}$ is Cauchy, and hence converges to some $b_m.$
Let
$$
	f(z) = \sum_m b_m z^m \in \F[[z,z^{-1}]].
$$
First, we need to check that $|f|_r<\infty,$ for $r_1\le r \le r_2.$
Since $a_{n,m}\to b_m,$ if $b_m\ne0,$ then
\hbox{$|a_{n,m}|=|b_m|$} for $n$ sufficiently large. Thus,
$$
	|b_m|r^m = |a_{n,m}|r^m \le |f_n|_r \le |f_n|_{\sup}
	\le \sup_n |f_n|_{\sup} < \infty.
$$
The last inequality follows from the assumption that $f_n$ is
Cauchy.
\smallskip
Second, we need to check that $|f_n-f|_{\sup}\to0.$ Because
$a_{n,m}\to b_m,$ for $n_m$ sufficiently large possibly depending on $m,$
we have
$$
	\sup_r |b_m-a_{n_m,m}|r^m < \varepsilon.
$$
On the other hand, since the $f_n$ are Cauchy, for $n'$ sufficiently large
and independent of $m,$ inequality~(\ref{cauchysupeqn}) is satisfied.
Therefore,
$$
	\sup_r\sup_m|b_m-a_{n',m}|r^m <
	\sup_r\sup_m \max\{|b_m-a_{n_m,m}|,|a_{n_m,m}-a_{n',m}|\}r^m
	< \varepsilon.\qedhere
$$
\end{proof}

\smallskip
By multiplicativity in Proposition~\ref{absrprop},  we can extend $|~|_r$ to
meromorphic functions.  We will need that meromorphic functions
are analytic away from poles.

\begin{prop}
\label{InvertProp}
Let $f$ be analytic on $\mathbf{B}_{\le r_1}$ with $r_1>0$ and $f(0)\ne0.$
Let
$$
        r_2 = \sup\{ r\le r_1 : |f|_r = |f(0)|\}>0.
$$
Then, there exists a unique analytic function
$g$ on $\mathbf{B}_{<r_2}$ such that $fg=1$ on $\mathbf{B}_{<r_2}.$
\end{prop}

\begin{proof}
The uniqueness of $g$ is clear.  By the non-Archimedean triangle inequality
and the choice of $r_2,$ we have 
$$
        \left | 1 - \frac{f}{f(0)}\right |_r < 1,
        \qquad\textrm{for all~} r<r_2.
$$
Hence,
$$
        \sum_{j=0}^\infty \left(1 - \frac{f}{f(0)}\right)^j
$$
converges to an analytic function $h$ on $\mathbf{B}_{<r_2}$ by
Proposition~\ref{CompleteProp}.  Here we are also using that a sequence
of functions converges in $\mathcal{A}(r_2)$
if and only if that sequence
converges in $\mathcal{A}[r]$ for all $r<r_2.$
Finally, set $g=h/f(0).$
\end{proof}

\subsubsection*{Liouville's Theorem Again and the Riemann Extension Theorem}

We now see how non-Archimedean analogs of Liouville's Theorem and
the Riemann Extension Theorem follow easily from 
the basic properties of $|~|_r.$  The proof of each proposition
is similar, so we state both propositions first and then give a joint
proof.

\begin{prop}[Liouville's Theorem]
\label{LiouvilleProp}
If $f$ is entire and $|f|_r$ is bounded for all
$r,$ then $f$ is constant.
\end{prop}

We say that an analytic function $f$ on $A[r_1,\infty)$ is \textbf{analytic
at infinity} if $f(1/z)$ is an analytic function on
$A[0,1/r_1]=\mathbf{B}_{\le r_1^{-1}}.$

\begin{prop}[Riemann Extension Theorem]
\label{RiemannProp}
If $f$ is analytic on $A[r_1,\infty)$ and the set
\hbox{$\{|f|_r : r\ge r_1\}$} is bounded in $\R,$
then $f$ is analytic at infinity.
\end{prop}

\begin{proof}[Proof of Propositions~\ref{LiouvilleProp} and
~\ref{RiemannProp}]
Write $f(z) = \sum c_n z^n.$  To prove either
proposition, we need to prove
that the existence of an index $n_0>0$ such that $c_{n_0}\ne0$
contradicts the boundedness of $|f|_r.$
But, if such a $n_0$ exits, then
\begin{equation*}
        |f|_r \ge |c_{n_0}| r^{n_0} \to \infty \textrm{~as~}
        r\to\infty. \qedhere
\end{equation*}
\end{proof}

An analytic function $f$ on $A[r_1,\infty)$ is said to be \textbf{meromorphic
at infinity} if $z^{-m}f(z)$ is analytic at infinity for some integer
$m\ge 0.$  If an analytic function $f$ on $A[r_1,\infty)$ is
not meromorphic at $\infty,$ then it is said to have 
\textbf{an essential singularity at infinity.}  We now state a proposition
that says if $|f|_r$ grows slowly as $r\to\infty,$
then it cannot have an essential singularity at infinity.

\begin{prop}
\label{meroinftyprop}
If $f$ is analytic on $A[r_1,\infty)$ and 
$$
        \limsup_{r\to\infty} \frac{\log |f|_r}{\log r} < \infty,
$$
then $f$ is meromorphic at infinity.
\end{prop}

\begin{proof}
Let $g(z)=z^{-m}f(z),$  which is also an analytic function on 
$A[r_1,\infty).$  Choose $m$ larger than
$$
        \limsup_{r\to\infty} \frac{\log |f|_r}{\log r}.
$$
Then,
$$
        \limsup_{r\to\infty} \log |g|_r < \infty.
$$
Hence $g$ is analytic at infinity by Proposition~\ref{RiemannProp}.
\end{proof}

\subsection{Valuation Polygons}

\subsubsection*{Preliminaries}

Let $f$ be analytic on $A[r_1,r_2].$  Write
$$
        f(z) = \sum_{n\in \Z} c_nz^n.
$$
For $r$ with $r_1\le r\le r_2,$ let
$$\fbox{$
        k(f,r) = \inf\{n\in\Z : |c_n|r^n = |f|_r\}\quad\textrm{and}\quad
        K(f,r) = \sup\{n\in\Z : |c_n|r^n = |f|_r\}.
$}$$
Note that in the special case of $r=0$ with $f(0)=0,$ we define
$$
	k(f,0) = 0 \qquad\textrm{and}\qquad
	K(f,0) = \inf\{n : c_n \ne 0\}.
$$
The integer $K(f,r)$ is often call the \textbf{central index}
and sometimes also has a role in complex analysis.

\smallskip
A radius $r$ such that $K(f,r)>k(f,r)$ is called a \textbf{critical
radius}.

\begin{prop}\label{CritRadDiscrete}
The set of critical radii for a Laurent series is discrete.
\end{prop}

\begin{proof}
Let $f$ be a Laurent series with critical radius $r',$
so $K(f,r')>k(f,r').$ Let $K=K(f,r')$ and $k=k(f,r').$
If $n> k,$ then either $a_n=0$ or
$|a_n|(r')^n\le|a_k|(r')^k.$ Hence, if $r<r',$ then
$$
	|a_n|r^n = \left(\frac{r}{r'}\right)^n|a_n|(r')^n
	\le \left(\frac{r}{r'}\right)^n|a_k|(r')^k
	=\left(\frac{r}{r'}\right)^{n-k}|a_k|r^k < |a_k|r^k,
$$
and so $K(f,r)\le k(f,r').$ Let $m$ be the largest integer $<k$
such that $a_m\ne0.$ If no such integer $m$ exists, then
$K(f,r)=k(f,r)=k$ for all $r<r',$ and so there are no critical
radii smaller than $r'.$ 
Otherwise, let $r''$ be the radius such that $|a_m|(r'')^m = |a_k|(r'')^k.$
Because,
$|a_m|(r')^m < |a_k|(r')^k,$ we know $r''<r'.$ Hence, for $r''<r<r',$
we have $K(f,r)=k(f,r),$ and so there are no critical radii
between $r''$ and $r'.$ By a similar argument, we see
that $k(f,r)\ge K(f,r')$ for $r\ge r',$ and if $M$ is the smallest
integer greater than $K$ such that $a_M\ne0,$ then $K(f,r)=k(f,r)=K$
for all $r'<r<r''',$ where $r'''>r'$ is such that
$|a_K|(r''')^K=|a_M|(r''')^M.$
\end{proof}

\begin{prop}\label{Kkmult}
If $f$ and $g$ are analytic on $A[r,r]$ then
$$
	K(fg,r)=K(f,r)+K(g,r)\qquad\textrm{and}\qquad
	k(fg,r)=k(f,r)+k(g,r).
$$
\end{prop}

\begin{proof}
We provide the proof for $K.$ The proof for $k$ is similar, or follows
by changing $z$ to $z^{-1}.$
Write $f(z)=\sum a_n z^n,$ $g(z)=\sum b_n z^n,$ and
$fg(z)=\sum c_n z^n.$ Let $m=K(f,r)+K(g,r).$ Then,
$$
	c_m = \sum_{i+j=m}a_ib_j.
$$
One of the terms in this sum comes from $i=K(f,r)$ and $j=K(g,r).$
In this case 
$$
	|a_ib_j|=\frac{|fg|_r}{r^m}.
$$
If $i<K(f,r),$ then $j>K(g,r),$
and so $|b_j|<|g|_r/r^j.$ But, $|a_i|\le|f|_r/r^i,$ and hence
$$
	|a_ib_j|<\frac{|fg|_r}{r^m}.
$$
Similarly, $|a_ib_j|<|fg|_r/r^m$ if $i>K(f,r).$
Hence, $|c_m|=|fg|_r/r^m.$ This shows that
$$
	K(fg,r)\ge K(f,r)+K(g,r).
$$
On the other hand, if $i+j>m,$ then either $i>K(f,r)$ or $j>K(g,r),$
and so we see \hbox{$K(fg,r)\le K(f,r) + K(g,r).$}
\end{proof}

\begin{prop}
\label{keqKprop}
Let $r>0.$  If $f$ is analytic on $A[r,r]$ with $K(f,r)=k(f,r),$ then
$f$ is invertible in $A[r,r].$
\end{prop}

\begin{proof}
Note $K(fz^m,r)=K(f,r)+m,$ and similarly for $k.$  Thus, by multiplying
$f$ by $z^{-k(f,r)},$ a unit in $A[r,r],$ we may assume
\hbox{$K(f,r)=k(f,r)=0.$}

Then, if $c_0$ is the constant term in the Laurent series defining $f,$
then the assumption that $K(f,r)=k(f,r)=0$ implies 
$$
	|f-c_0|_r < |c_0| \quad\textnormal{or in other words}\quad
	|c_0^{-1}f - 1|_r < 1.
$$
Thus,
\begin{equation*}
	f^{-1} = c_0^{-1}[1-(1-c_0^{-1}f)]^{-1}
	= c_0^{-1}[1+(1-c_0^{-1}f)^2+(1-c_0^{-1}f)^3+\dots] \qedhere
\end{equation*}
\end{proof}

\subsubsection*{Valuation Polygons for Polynomials}

Consider a monic linear polynomial $L(z)=z-a.$
Then, clearly
$$
	|L|_r = \left\{\begin{array}{ll}
	|a| &\textrm{~if~}r\le|a|\\
	r &\textrm{~if~}r\ge|a|.
	\end{array}\right.
$$
The following figure is a graph of $\log|L|_r$ as a function of $\log r.$

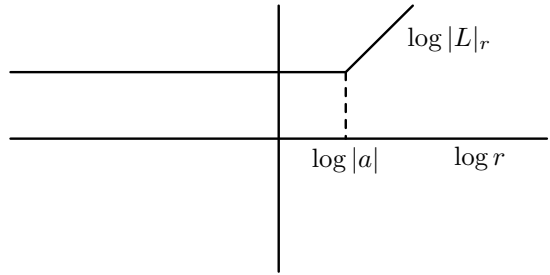
\smallskip
\centerline{\epsfig{file=pfunc1.ps}}

\smallskip\noindent
Notice that the corner of the graph indicates that $L(z)$ has a zero
with $|z|=|a|.$

\smallskip
Now suppose that $P(z)=(z-a)^n(z-b)^m$ with $0<|a|<|b|.$
Then,
$$
	\log|P|_r = \left\{\begin{array}{cc}
	n\log|a|+m\log |b|&\textrm{~if~}r\le|a|\\
	n\log r + m\log |b| &\textrm{~if~}|a|\le r \le |b|\\
	(n+m)\log r &\textrm{~if~}r\ge|b|.
	\end{array}\right.
$$
This time the graph of $\log|P|_r$ as a function of $\log r$
looks like

\smallskip
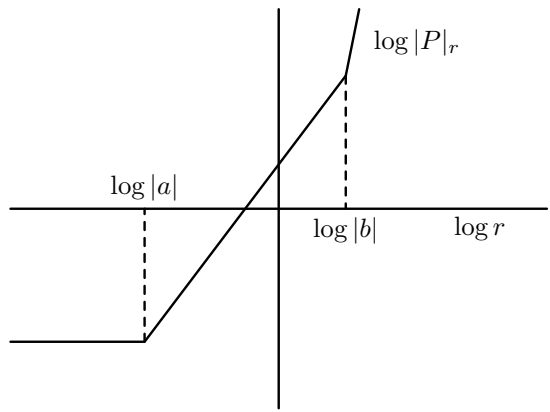
\centerline{\epsfig{file=pfunc2.ps}}

\smallskip\noindent
Again, we see a piecewise linear graph whose corners indicate the
location of the zeros of $P.$ Notice also that the change in slope
indicates the number of zeros with that absolute value: in this
case the slope goes from $0$ to $n$ at $r=|a|$ and from $n$ to $n+m$
at $r=|b|.$ The above graph is called the \textbf{valuation polygon}
of $P.$  Observe that for the above example,
$K(P,|a|)=n$ and $k(P,|a|)=0,$ that $K(P,|b|)=m+n$ and
$k(P,|b|)=n,$ and that $K(P,r)=k(P,r)$ for all $r\ne|a|,|b|.$
Thus, the corners of the valuation polygon correspond to the
critical radii.

\smallskip
If $P$ is an arbitrary polynomial, we can write
$$
	P(z)=cz^{m_0}\prod_j(z-a_j)^{m_j},
$$
and we see that $\log r \mapsto \log|P|_r$ is a piecewise linear
function whose corners indicate the locations of the zeros at $P,$
and such that the change in slope at the corners indicates the number
of zeros, counting multiplicity, that $P$ has at that absolute value.

\begin{remark*} What is known as the \textbf{Newton polygon} is a
polygon dual to the valuation polygon in a certain sense. As such,
the Newton polygon of a polynomial $P$ also encodes the locations of the
zeros of $P,$ but I prefer the valuation polygon to the Newton polygon
for non-Archimedean function theory.  See \cite[pp.~300]{Robert}
for a more detailed
description of the Newton polygon and its relationship to the valuation 
polygon.
\end{remark*}

We now show that for a polynomial $P,$
the corners of the valuation polygon occur precisely
at the critical radii and that $P$ has $K(P,r)-k(P,r)$
zeros, counting multiplicity, with absolute value $r.$

\begin{prop}\label{PolyCrit}
A polynomial $P$ has $K(P,r)-k(P,r)$ zeros, counting multiplicity,
of absolute value $r.$
\end{prop}

\begin{proof}
Without loss of generality, we may assume $P$ is monic,
and we proceed by induction.  
The proof essentially amounts to Gauss's Lemma.  

\smallskip
The proposition
is clear if $P$ is a linear polynomial.  Now consider
$$
	P(z) = z^n + b_{n-1}z^{n-1}+\dots+b_0 = 
	(z-z_1)(z^{n-1}+a_{n-2}z^{n-2}+\dots+a_0).
$$
Let
$$
        k_1 = \inf\{m : |a_m|r^m = \sup_\ell |a_\ell| r^\ell\}
	\qquad\textrm{and}\qquad
        K_1 = \sup\{m : |a_m|r^m = \sup_\ell |a_\ell| r^\ell\}.
$$
By the induction hypotheses, $z^{n-1}+a_{n-2}z^{n-2}+\dots+a_0$
has $K_1-k_1$ zeros on $|z|=r.$

\smallskip
Now, $|b_0|=|a_0||z_1|.$  Hence,
\begin{equation}\label{b0k1pos}
	|b_0| < \frac{|z_1|}{r}|a_{K_1}|r^{K_1+1} \qquad\textnormal{if~}k_1>0,
\end{equation}
and
\begin{equation}\label{b0k1zer}
	|b_0| = \frac{|z_1|}{r}|a_{K_1}|r^{K_1+1} \qquad\textnormal{if~}k_1=0.
\end{equation}

For $0<j<k_1$ and for $K_1+1<j\le n,$
\begin{eqnarray}\label{bjltk1}
	|b_j|r^j &=& |a_{j-1}-z_1a_j|r^j\nonumber\\
	&\le& \max\left\{|a_{j-1}|r^{j-1},\frac{|z_1|}{r}|a_j|r^j\right\}r
	\nonumber\\
	&<& \max\left\{1,\frac{|z_1|}{r}\right\}|a_{K_1}|r^{K_1+1}.
\end{eqnarray}

\smallskip
If $j=k_1>0$ and $|z_1|<r,$ then
\begin{eqnarray}\label{bk1z1ltr}
	|b_{k_1}|r^{k_1} &=& |a_{k_1-1}-z_1a_{k_1}|r^{k_1}\nonumber\\
	&\le& \max\left\{|a_{k_1-1}|r^{k_1-1},
	\frac{|z_1|}{r}|a_{k_1}|r^{k_1}\right\}r\nonumber\\
	&<& |a_{K_1}|r^{K_1+1}.
\end{eqnarray}
Similarly, if $j=k_1>0$ and $|z_1|\ge r,$ then
\begin{eqnarray}\label{bk1z1ger}
	|b_{k_1}|r^{k_1} &=& |a_{k_1-1}-z_1a_{k_1}|r^{k_1}\nonumber\\
	&=& \max\left\{|a_{k_1-1}|r^{k_1-1},
	\frac{|z_1|}{r}|a_{k_1}|r^{k_1}\right\}r\nonumber\\
	&=& \frac{|z_1|}{r}|a_{K_1}|r^{K_1+1}.
\end{eqnarray}

\smallskip
For $k_1+1\le j \le K_1,$ 
\begin{eqnarray}\label{bjk1plus1toK1}
	|b_j|r^j &=& |a_{j-1}-z_1a_j|r^j\nonumber\\
	&\le& \max\left\{|a_{j-1}|r^{j-1},\frac{|z_1|}{r}|a_j|r^j\right\}r
	\nonumber\\
	&\le& \max\left\{1,\frac{|z_1|}{r}\right\}|a_{K_1}|r^{K_1+1}.
\end{eqnarray}
For $j=k_1+1$ and $|z_1|<r,$ then
\begin{eqnarray}\label{bk1plus1z1ltr}
	|b_{k_1+1}|r^{k_1+1} &=& |a_{k_1}-z_1a_{k_1+1}|r^{k_1+1}\nonumber\\
	&=& \max\left\{|a_{k_1}|r^{k_1},
		\frac{|z_1|}{r}|a_{k_1+1}|r^{k_1+1}\right\}r
	\nonumber\\
	&=& |a_{K_1}|r^{K_1+1}.
\end{eqnarray}
For $j=K_1$ and $|z_1|>r,$ then
\begin{eqnarray}\label{bK1z1gtr}
	|b_{K_1}|r^{K_1} &=& |a_{K_1-1}-z_1a_{K_1}|r^{K_1}\nonumber\\
	&=& \max\left\{|a_{K_1-1}|r^{K_1-1},
		\frac{|z_1|}{r}|a_{K_1}|r^{K_1}\right\}r
	\nonumber\\
	&=& \frac{|z_1|}{r}|a_{K_1}|r^{K_1+1}.
\end{eqnarray}

\smallskip
For $j=K_1+1<n$ and $|z_1|> r,$ then
\begin{eqnarray}\label{bK1pl1z1gtr}
	|b_{K_1+1}|r^{K_1+1} &=& |a_{K_1}-z_1a_{K_1+1}|r^{K_1+1}\nonumber\\
	&\le& \max\left\{|a_{K_1}|r^{K_1},\frac{|z_1|}{r}|a_{K_1+1}|
	r^{K_1+1}\right\}r\nonumber\\
	&<& \frac{|z_1|}{r}|a_{K_1}|r^{K_1+1}. 
\end{eqnarray}
Similarly, for $j=K_1+1<n$ and $|z_1| \le r,$ 
\begin{eqnarray}\label{bK1pl1z1ler}
	|b_{K_1+1}|r^{K_1+1} &=& |a_{K_1}-z_1a_{K_1+1}|r^{K_1+1}\nonumber\\
	&=& \max\left\{|a_{K_1}|r^{K_1},\frac{|z_1|}{r}|a_{K_1+1}|
	r^{K_1+1}\right\}r\nonumber\\
	&=& |a_{K_1}|r^{K_1+1}. 
\end{eqnarray}

\smallskip
If $K_1=n-1,$ then
\begin{equation}\label{K1max}
	|b_n|r^n = r^n = |a_{K_1}|r^{K_1+1}.
\end{equation}

\smallskip
If $|z_1|=r,$ then by $(\ref{b0k1pos})$ and $(\ref{bjltk1}),$
we have
$$
	|b_j|r^j < |a_{K_1}|r^{K_1+1},\qquad
	\textnormal{for~}0\le j < k_1 \textnormal{~and~for~}
	K_1+1<j\le n.
$$
Also, by $(\ref{b0k1zer}),$ $(\ref{bk1z1ger}),$ $(\ref{bK1pl1z1ler}),$ 
and $(\ref{K1max}),$ we have
$$
	|b_{k_1}|r^{k_1} \;=\; |b_{K_1+1}|r^{K_1+1} \;=\; |a_{K_1}|r^{K_1+1}.
$$
Hence, taking $(\ref{bjk1plus1toK1})$ into account,
$K=K_1+1,$ while $k=k_1,$ and so 
$$
	K-k=K_1-k_1+1,
$$ 
as was to be shown.

\smallskip
If $|z_1| < r,$ then by $(\ref{b0k1pos})$ and $(\ref{bjltk1}),$
we have
$$
	|b_j|r^j < |a_{K_1}|r^{K_1+1},\qquad
	\textnormal{for~}0\le j < k_1 \textnormal{~and~for~}
	K_1+1<j\le n.
$$
Also, by $(\ref{b0k1zer})$ and $(\ref{bk1z1ltr}),$ 
we have
$$
	|b_{k_1}|r^{k_1} \;<\; |a_{K_1}|r^{K_1+1}.
$$
Finally, by $(\ref{bK1pl1z1ler}),$ $(\ref{K1max}),$
and $(\ref{bk1plus1z1ltr}),$ we have
$$
	|b_{K_1+1}|r^{K_1+1} \;=\; |a_{K_1}|r^{K_1+1}
	\;=\; |b_{k_1+1}|r^{k_1+1}.
$$
Thus, again taking $(\ref{bjk1plus1toK1})$ into account,
$K=K_1+1,$ and $k=k_1+1,$ so $K-k=K_1-k_1$ as was to be shown.

\smallskip
For the last case in the induction, suppose $|z_1|>r.$
By $(\ref{b0k1pos}),$ $(\ref{bjltk1}),$ $(\ref{bK1pl1z1gtr}),$
and $(\ref{K1max}),$ we have
$$
	|b_j|r^{j} < \frac{|z_1|}{r}|a_{K_1}|r^{K_1+1}\qquad
	\textnormal{for~}j<k_1 \textnormal{~and~for~}j>K_1.
$$
By $(\ref{bjk1plus1toK1}),$
$$
	|b_j|r^{j} \le \frac{|z_1|}{r}|a_{K_1}|r^{K_1+1}\qquad
	\textnormal{for~}k_1+1 \le j \le K_1-1.
$$
By $(\ref{b0k1zer}),$ $(\ref{bk1z1ger}),$ and $(\ref{bK1z1gtr}),$ we have
$$
	|b_{k_1}|r^{k_1} = |b_{K_1}|r^{K_1} = 
	\frac{|z_1|}{r}|a_{K_1}|r^{K_1+1}.
$$
Thus, $k=k_1$ and $K=K_1,$ so $K-k=K_1-k_1$ as required.
\end{proof}

\subsection{Euclidean Division Algorithm}

Following \cite{Amice}, we analyze the Euclidean division algorithm
for Laurent series.

\smallskip
Let $r\ge0.$  A polynomial $P$ is called \textbf{$r$-dominant} if
$K(P,r)=\deg P$ and it is called \textbf{$r$-extremal} if it is
$r$-dominant and in addition $k(P,r)=0.$ Thus, a polynomial is 
$r$-dominant if and only if all of its zeros are located in
$\mathbf{B}_{\le r},$ and a polynomial is $r$-extremal if and only if
all of its zeros are located in the annulus $|z|=r.$

\begin{lemma}[Continuity of Division]\label{EuclidCont}
Let $r>0$ and let $f$ be analytic on $\mathbf{B}_{\le r}.$
Let $P$ be a polynomial in $\mathbf{F}[z]$ with $P\not\equiv0,$
let $r>0,$ and
assume that $P$ is $r$-dominant. Then there exist a unique
function $q$ analytic on $\mathbf{B}_{\le r}$
and a unique polynomial $R$ such that
\begin{quote}\begin{enumerate}
\item[\textup{(i)}] $f=Pq+R;$
\item[\textup{(ii)}] $\deg R < \deg P;$
\item[\textup{(iii)}] $|R|_r \le |f|_r;$ and
\item[\textup{(iv)}] $\ds |q|_r \le \frac{|f|_r}{|P|_r}.$
\end{enumerate}\end{quote}
\end{lemma}

\begin{remark*} Lemma~\ref{EuclidCont} is referred to as 
continuity of division because it implies that if $f_1$ and $f_2$
are analytic functions with $|f_1-f_2|_r$ small, and if $f_1$ and $f_2$
are each divided by $P$ to get quotients $q_1$ and $q_2$ and remainders
$R_1$ and $R_2,$ then $|q_1-q_2|$ and $|R_1-R_2|$ are both small.
\end{remark*}

\begin{proof}
We first prove the case when $f$ is also a polynomial.
In that case, the Euclidean algorithm gives unique
polynomials $q$ and $R$ satisfying (i) and (ii),
so it remains to check (iii) and (iv).

\smallskip
First consider the special case that $r=1$ and all the coefficients of
$f$ and $P$ have absolute value at most $1,$ 
\textit{i.e.,} are elements of $\O,$ and that at least one coefficient
in each polynomial has absolute value $1.$ This means
$|f|_1=|P|_1=1,$ and since $P$ is $1$-dominant, its leading coefficient
must have absolute value $1.$ Then, the Euclidean division algorithm
produces polynomials $R$ and $q$ with coefficients in $\O$ and since 
$|f|_1=|P|_1=1,$ (iii) and (iv) follow. Since multiplying $f$ by
a constant multiplies $q$ and $R$ by the same constant, the result continues
to hold without the assumption that the coefficients of $f$
are in $\O.$ Multiplying $P$
by a constant divides $q$ by the same constant and does not change $R,$
so the result continues to hold without the 
assumption that the coefficients of $P$
are in $\O$ as well.

\smallskip
Still in the case that $f$ is a polynomial, 
if $r$ is in $|\Ftimes|,$ then by choosing $a$ in $\F$ with $|a|=r$
and changing variables by replacing $z$ with $az,$ we reduce to the case
above. If $|~|$ is trivial, then the lemma is also trivial.
When $|~|$ is non-trivial, if $r$ is not in $|\Ftimes|,$
then for $r'$ in $|\Ftimes|$ 
sufficiently close
to $r,$ we will have that $P$ is $r'$-dominant, and the lemma follows
since $|~|_r$ is continuous in $r.$

\smallskip
If $f$ is not a polynomial, we can find a sequence of polynomials
$f_n$ such that \hbox{$|f-f_n|_r\to0,$} for instance by truncating the power
series representation of $f$ to higher and higher orders.
Letting $q_n$ and $R_n$ be the quotients and remainders obtained by
dividing the $f_n$ by $P,$ we have by the polynomial version of the
lemma already proven that $q_n$ and $R_n$ are sequences
in $\mathcal{A}[r]$ that are Cauchy sequences with respect to $|~|_r.$
Therefore they converge to $q$ and $R$ in $\mathcal{A}[r]$
by Proposition~\ref{CompleteProp}.
As $\deg R_n < \deg P,$ the $R_n$ must converge to a polynomial, also of
degree $<\deg P.$ Properties~(i)--(iv) are preserved under taking limits
as $n\to\infty.$ Property~(iv) ensures that the quotient $q$ is analytic
on $\mathbf{B}_{\le r}.$

\smallskip
Now to check uniqueness in the general case, suppose we have
$$
	Pq_1+R_1 = f = Pq_2+R_2.
$$
Then, we would have \hbox{$P(q_1-q_2)=R_2-R_1.$}
Hence, if $q_1\ne q_2,$ then
$$
K(R_2-R_1,r)=
K(P(q_1-q_2),r)=K(P,r)+K(q_1-q_2,r)=\deg P + K(q_1-q_2,r)\ge \deg P,
$$
which contradicts the fact that $R_2-R_1$ is a polynomial of degree less
than $\deg P.$
\end{proof}

\begin{cor}\label{LaurentDivision} Let $r_1\le r \le r_2$ with $r>0,$
let $f$ be in $\mathcal{A}[r_1,r_2],$ and let $P$ be an $r$-extremal
polynomial. Then, there exists a unique
Laurent series $q$ in $\mathcal{A}[r_1,r_2]$ and a unique polynomial $R$
such that
\begin{quote}\begin{enumerate}
\item[\textup{(i)}]$f=Pq+R;$
\item[\textup{(ii)}]$\deg R < \deg P;$
\item[\textup{(iii)}] $|R|_r \le |f|_r;$ and
\item[\textup{(iv)}] $\ds |q|_r \le \frac{|f|_r}{|P|_r}.$
\end{enumerate}\end{quote}
\end{cor}

\begin{proof}
We begin by proving uniqueness. Suppose
$$
	Pq+R = P\tilde q + \widetilde{R}.
$$
Then,
$$
	P(q-\tilde q) = \widetilde{R}-R.
$$
Hence,
\begin{align*}
	K(\widetilde{R}-R,r) &= K(P,r)+K(q-\tilde q,r) 
		= \deg P + K(q-\tilde q,r) \qquad\textrm{and}\\
	k(\widetilde{R}-R,r) &= k(P,r)+k(q-\tilde q,r) 
		= 0 + k(q-\tilde q,r)
\end{align*}
since $P$ is $r$-extremal. Hence,
$$
	K(\widetilde{R}-R,r)-k(\widetilde{R}-R,r)\ge \deg P,
$$
which contradicts the fact that $\widetilde{R}-R$ is a polynomial
of degree $<\deg P.$

\smallskip
To prove existence, write
$$
	f(z) = \sum_{n=0}^\infty a_n z^n + \sum_{n=-\infty}^{-1} a_n z^n
	= f_+(z) + f_-(z).
$$
First, apply the division algorithm to $f_+$ to get a power series
$q_+$ and a polynomial $R_+$ such that
\hbox{$f_+=Pq_++R_+.$} From the lemma, we know $\deg R_+ < \deg P,$
$$
	|q_+|_r \le \frac{|f_+|_r}{|P|_r} \qquad\textrm{and}\qquad
	|R_+|_r \le |f_+|_r.
$$
Since $r_2\ge r,$ $P$ is also $r_2$ dominant, and the lemma also implies that
$q_+$ is in $\mathcal{A}[r_2].$
Next, observe that $z^{\deg P}P(z^{-1}),$ the palindrome of $P,$
is $r^{-1}$-extremal. Because $f_-(z^{-1})/z$ is a power series
that converges for $|z|=r^{-1},$ we can apply the division algorithm to get
a power series $q_-$ and a polynomial $R_-$ with $\deg R_- < \deg P$
such that
$$
	z^{\deg P}\frac{f_-(z^{-1})}{z}=z^{\deg P}P(z^{-1})q_-(z)+R_-(z).
$$
Hence,
$$
	f_-(z)=P(z)\frac{q_-(z^{-1})}{z}+z^{\deg P-1}R_-(z^{-1}).
$$
From the lemma, we also know that
$$
	|z^{\deg P-1}R_-(z^{-1})|_r = |z^{1-\deg P}R_-(z)|_{r^{-1}}
	\le|f_-(z^{-1})|_{r^{-1}} = |f_-(z)|_r,
$$
and
$$
	|z^{-1}q_-(z^{-1})|_r = |zq_-(z)|_{r^{-1}} \le
	\frac{|f_-(z^{-1})|_{r^{-1}}}{|P(z^{-1})|_{r^{-1}}}
	=\frac{|f_-(z)|_r}{|P(z)|_r}.
$$
Because $r_1\le r,$  $z^{\deg P}P(z^{-1})$ is $r_1^{-1}$ dominant,
and the lemma implies $z^{-1}q_-(z^{-1})$ is in 
$\mathcal{A}[r_1,\infty).$
Thus, the corollary follows by setting
$$
	q(z)=q_+(z)+z^{-1}q_-(z^{-1}) \qquad\textrm{and}\qquad
	R(z)=R_+(z)+z^{\deg P-1}R_-(z^{-1}). \qedhere
$$
\end{proof}

\subsubsection*{Weierstrass Preparation}

The following theorem is a non-Archimedean one variable version
of the Weierstrass Preparation Theorem and is sometimes referred
to as a version of Hensel's Lemma.

\begin{theorem}[Weierstrass Preparation]
\label{WeierstrassThm}
Let $f$ be an analytic function on $A[r_1,r_2].$  Let $r$ be such
that \hbox{$r_1\le r \le r_2.$}  Let \hbox{$d=K(f,r)-k(f,r).$}
Then, there exists a unique pair $(P,u)$ such that $f=Pu,$ such that
$P$ is a polynomial
of degree $d$ with $P(0)=1,$ $k(P,r)=0,$ and $K(P,r)=d,$ and such that
$u$ is analytic on $A[r_1,r_2]$ with $k(u,r)=K(u,r).$
\end{theorem}

\begin{proof}
Multiplying $f$ by a constant and a suitable power of $z,$ we may,
without loss of generality, assume that $|f|_r=1$ and $k(f,r)=0.$
We show existence by an inductive construction.
Write
$$
	f(z)=\sum_{n=-\infty}^\infty a_n z^n \qquad\textrm{and let}
	\qquad P_1(z) = \sum_{n=0}^d a_n z^n,
$$
where $d=K(f,r).$ Note that $P_1$ is $r$-extremal, and let
$q_1$ and $R_1$ be the quotient and remainder
provided by Corollary~\ref{LaurentDivision}. Observe that
$$
	f-P_1 = P_1(q_1 - 1) + R_1,
$$
and by the uniqueness of the division algorithm and
Corollary~\ref{LaurentDivision}, we have
$$
	|R_1|_r \le |f-P_1|_r < 1,
$$
keeping in mind that $|P_1|_r=1.$

\smallskip
Now assume that 
for $i=1,\dots,n$ we 
have found $(P_i,q_i,R_i)$ where the $P_i$ are degree $d$
$r$-extremal polynomials with $|P_i|_r=1,$
where the $R_i$ are polynomials of degree $<d,$
where the $q_i$ are analytic on $A[r_1,r_2]$ with $|q_i|_r=1,$
where $f=P_iq_i+R_i,$ and where the following
inequalities hold:
\begin{enumerate}
\item[(a)] $|R_i|_r \le |f-P_1|^i$ for $i=1,\dots,n;$
\item[(b)] $|P_i-P_{i-1}|_r \le |f-P_1|^{i-1}$ for $i=2,\dots,n;$ and
\item[(c)] $|q_i-q_{i-1}|_r \le |f-P_1|^i$ for $i=2,\dots,n.$
\end{enumerate}
Set $P_{n+1}=P_n+R_n.$ Now, $P_{n+1}$ and $P_n$ have the same
top degree term and \hbox{$|R_n|_r < 1 = |P_n|_r,$}
so $|P_{n+1}|_r=1$ and $P_{n+1}$ is $r$-dominant.
Also, 
$$
	|R_n(0)| \le |R_n|_r < 1 = |P_n|_r = |P_n(0)|,
$$
and so by Exercise~\ref{sumeqmax}, $k(P_{n+1},r)=0,$
and $P_{n+1}$ is $r$-extremal. Let $q_{n+1}$ and $R_{n+1}$
be the quotient and remainder obtained by
dividing $f$ by $P_{n+1}.$ Now,
$$
	|P_{n+1}-P_n|_r = |R_n|_r \le |f-P_1|_r^n,
$$
and so (b) is satisfied by $P_{n+1}.$ 

\smallskip
Re-arranging the equation
$$
	P_nq_n+R_n = f = P_{n+1}q_{n+1}+R_{n+1} = (P_n+R_n)q_{n+1}+R_{n+1},
$$
we get
\begin{align}
\label{wpeqone}
	-R_nq_{n+1}&=P_n(q_{n+1}-q_n)+R_{n+1}-R_n \\
\label{wpeqtwo}
	\textrm{and}\qquad
	R_n(1-q_{n+1})&=P_n(q_{n+1}-q_n)+R_{n+1}.
\end{align}
Applying Corollary~\ref{LaurentDivision} to $(\ref{wpeqone}),$ we have
$$
	|q_{n+1}-q_n|_r \le |R_n|_r\cdot|q_{n+1}|_r\le|f-P_1|_r^n\cdot 1<1,
$$
and hence $|q_{n+1}|_r=|q_n|_r=1$ by Exercise~\ref{sumeqmax}.
Now,
$$
	|1-q_{n+1}|_r = |1-q_1+q_1-q_2+\dots+q_n-q_{n+1}|_r
	\le \max\{|1-q_1|_r,|q_1-q_2|_r,\dots,|q_n-q_{n+1}|_r\}
	\le |f-P_1|_r.
$$
Combining this with applying Corollary~\ref{LaurentDivision}
to $(\ref{wpeqtwo}),$ we get
$$
	|q_{n+1}-q_n|_r \le |R_n|_r\cdot|1-q_n|_r\le
	|f-P_1|_r^n\cdot|f-P_1|_r = |f-P_1|_r^{n+1},
$$
which shows~(c), and~(a) follows similarly.

\smallskip
By~(a), $\lim R_n=0.$ Let $P=\lim P_n$ and $u=\lim q_n.$
Then, $f=Pu$ and $P$ is an $r$-extremal degree $d$ polynomial
since each $P_n$ is.
As $u$ is the quotient under long division, Corollary~\ref{LaurentDivision}
implies that $u$ is in $\mathcal{A}[r_1,r_2].$ Because
$$
	d=K(f,r)-k(f,r)=K(P,r)-k(P,r)+K(u,r)-k(u,r)=d+K(u,r)-k(u,r),
$$
we see that $K(u,r)=k(u,r),$ as was to be shown.

\smallskip
Only the uniqueness remains. Suppose $Pu=\widetilde{P}\tilde u.$
By Proposition~\ref{keqKprop}, $u$ is invertible in $\mathcal{A}[r,r].$
Hence,
$P=u^{-1}\widetilde{P}\tilde u.$ However, $u^{-1}\tilde u$ is the
quotient of $P$ by $\widetilde{P},$ and is hence a polynomial, by
the uniqueness of the quotient. Since $P$ and $\widetilde{P}$ have the
same degree, this implies $u^{-1}\tilde u$ is constant. Since
$P(0)=\widetilde{P}(0)=1,$ that constant is $1$ and we conclude
$u=\tilde u$ and $P=\widetilde{P}.$
\end{proof}

Theorem~\ref{WeierstrassThm} allows us to connect $|f|_r,$
$K(f,r),$ $k(f,r),$  and the locations of the zeros of $f.$

\begin{theorem}
\label{NewtonPolyThm}
Let $f$ be analytic on $A[r_1,r_2].$
If $r_1\le\rho\le R\le r_2,$
then $f$ has 
$$K(f,R)-k(f,\rho)$$ zeros in 
$A[\rho,R]$ counting multiplicity.
\end{theorem}

\begin{proof}
Clearly, $K(f,r)$ and $k(f,r)$ are non-decreasing in $r,$ so it suffices
to prove the theorem for $\rho=R=r.$  
Also, the case of $r=0$ is clear,
so we assume $r>0.$

\smallskip
Write $f=Pu$ as in Theorem~\ref{WeierstrassThm}.
From Proposition~\ref{PolyCrit}, $P$ has 
$$K(P,r)-k(P,r)=K(f,r)-k(f,r)
$$
zeros, all with absolute value $r.$
\end{proof}

Theorem~\ref{NewtonPolyThm} tells us that the valuation polygon
for a Laurent series $f$ works just like the polynomial case.
Namely, the zeros of $f$ occur precisely at the corners of the valuation
polygon, and the sharpness of the corner determines the number of
zeros, counting multiplicity.

\smallskip
The connection between the locations of the zeros of a non-Archimedean
analytic function and its Laurent series coefficients given by 
Theorem~\ref{NewtonPolyThm} is strikingly different from the classical
complex case and is responsible for most of the differences between 
classical function theory and non-Archimedean function theory.

\begin{cor}[Identity Principle]
\label{IdentityCor}
If $f\not\equiv0$ is analytic on $A[r_1,r_2],$ with
\hbox{$r_2<\infty,$} then $f$ has at most finitely
many zeros in $A[r_1,r_2].$
\end{cor}

\begin{proof}
The numbers $k(f,r)$ and $K(f,r)$  are non-decreasing functions of $r,$
so the number of zeros is bounded by $K(f,r_2)-k(f,r_1).$
\end{proof}

Theorem~\ref{NewtonPolyThm} allows us to immediately conclude
non-Archimedean analogs of Picard's theorems for maps to the projective
line. 

\begin{cor}[Little Picard]
\label{LittlePicCor}
If $f$ is analytic and zero free
on $A[0,\infty),$ then $f$ is constant.
\end{cor}

\begin{proof} Since $f$ is zero free, $f(0)\ne0,$ and so
$$
        f(z)=a_0+\sum_{n=1}^\infty a_nz^n, \quad\textrm{with~}a_0\ne0.
$$
Since we don't have any zeros, we have by Theorem~\ref{NewtonPolyThm},
that
$$
        \sup_{n\ge1}|a_n|r^n < |a_0| \quad\textrm{for all~}r.
$$
This is clearly impossible unless $a_n=0$ for all $n\ge1.$
\end{proof}

\begin{cor}[Big Picard]
\label{BigPicCor}
If $f$ is analytic and zero-free on $A[r_1,\infty),$ then
$$
        \limsup_{r\to\infty} \frac{\log |f|_r}{\log r} < \infty.
$$
In particular, $f$ is meromorphic at infinity.
\end{cor}

\begin{remark*}
Corollary~\ref{BigPicCor} seems to have first appeared in the literature
in \cite{vanderPut}, although van~der~Put himself says that this surely
must have been known much earlier.
\end{remark*}

\begin{proof}
Write 
$$
        f(z) = \sum_{n\in\Z}a_nz^n.
$$
Because $f$ is zero free,
$$
        |a_n|r^n < |f|_{r_1} \quad \textrm{for~}|n| \textrm{~large.}
$$
Thus, $f$ clearly has only finitely many non-zero $a_n$ with $n>0.$
The proof is completed by Proposition~\ref{meroinftyprop}.
\end{proof}

\subsection{Poisson-Jensen Formula}
We conclude this lecture with a non-Archimedean Poisson-Jensen
formula. Let $f$ be an analytic function on $A[r_1,r_2]$
which is not identically zero.
We define the \textbf{counting function} $N(f,0,r)$ by defining
$$
        N(f,0,r) = \!\!\!\!\!\!\sum_{\begin{array}{c}\scriptstyle0\ne z\in A[r_1,r]\\
        \scriptstyle\textrm{s.t.~}f(z)=0\end{array}} \!\!\!\!\!\!
	\log \frac{r}{|z|}.
$$
Here, we count the zeros of $f(z)$ with multiplicity. 
If $r_1=0,$ it is convenient to add the term $K(f,0)\log r$
to the definition of $N(f,0,r).$
By the identity principle, the sum defining
$N$ is finite if $r\in[r_1,r_2].$  \textit{Note that $N$ implicitly depends
on the lower radius in the annulus $r_1.$}

\begin{theorem}[Poisson-Jensen]
\label{PoissonJensen}
Let $f$ be a non-constant analytic function on $A[r_1,r_2),$ with
$r_2\le\infty.$  Let
$$
	f(z)=\sum_{n\in\Z}a_nz^n	
$$
be the Laurent expansion for $f.$
Then, for all $r\in[r_1,r_2),$ we have
$$
        N(f,0,r) +k(f,r_1)\log r + \log|a_{k(f,r_1)}|
	=   \log |f|_r\qquad\textnormal{if~}
	r_1>0
$$
or
$$
	N(f,0,r)+\log|a_{K(f,0)}| = \log |f|_r \qquad\textnormal{if~}
	r_1=0.
$$
\end{theorem}

\begin{remark*}
If $r_1=0$ or $r_2<\infty,$ then the theorem says the difference between
$N(f,0,r)$ and $\log |f|_r$ remains bounded as $r\to r_2.$  
If $r_1>0$ and $r_2=\infty,$ then the difference is bounded by
$O(\log r)$ as $r\to\infty.$
\end{remark*}

\begin{proof}
This is basically unwinding definitions and understanding the valuation
polygon.
Recall that by Proposition~\ref{CritRadDiscrete},
there are only finitely many critical
points in $[r_1,r]$ if $r<r_2.$

\smallskip
In the case $r_1=0,$ let $r'$ be the smallest positive critical point.
By Theorem~\ref{WeierstrassThm}, $f$ has no zeros with absolute value
between $0$ and $r'.$  Thus, for $0<r\le r',$
\begin{eqnarray*}
	N(f,0,r) \;=\; K(f,0)\log r&=&K(f,r)\log r \\
	&=& \log|f|_r-\log|a_{K(f,r)}|
	\;=\; \log|f|_r - \log|a_{K(f,0)}|,
\end{eqnarray*}
where we have used $K(f,r)=K(f,0).$

\smallskip
In the case $r_1>0,$ let $r'$ be the smallest critical point $>r_1.$
Again, by Theorem~\ref{WeierstrassThm}, $f$ has no zeros with absolute
value between $r_1$ and $r'.$  Therefore, for \hbox{$r_1\le r \le r',$}
again using $K(f,r)=K(f,r_1),$ and the fact that
\begin{equation}\label{Kkr1}
	\log|a_{k(f,r_1)}| + k(f,r_1)\log r_1 = \log |f|_{r_1} = 
	\log|a_{K(f,r_1)}| + K(f,r_1)\log r_1,
\end{equation}
we get
\begin{eqnarray*}
	N(f,0,r) &=& \!\!\!\!\!\!\!\!
	\sum_{\begin{array}{c}\scriptstyle |z|=r_1\\
        \scriptstyle\textrm{s.t.~}f(z)=0\end{array}} 
	\!\!\!\!\!\!\!\!\log \frac{r}{|z|}\\
	&=&[K(f,r_1)-k(f,r_1)]\log\frac{r}{r_1}\\
	&=&K(f,r_1)\log r -k(f,r_1)\log r + 
		k(f,r_1)\log r_1 -K(f,r_1)\log r_1 \\
  \textrm{[From~}(\ref{Kkr1})\textrm{]}
	&=&K(f,r)\log r - k(f,r_1)\log r +\log|a_{K(f,r_1)}|
	   -\log|a_{k(f,r_1)}|\\
	&=&\log|f|_r-\log|a_{K(f,r)}| -k(f,r_1)\log r \\
	&&\qquad+\log|a_{K(f,r_1)}|
	   -\log|a_{k(f,r_1)}|\\
	&=& \log|f|_r-k(f,r_1)\log r-\log|a_{k(f,r_1)}|.
\end{eqnarray*}

\smallskip
Thus, in both cases, we see that the desired formula is correct 
for $r$ between $r_1$ and the first critical point bigger than $r_1.$
We then simply need to check that we can pass through each critical point.
Thus, assume $r'$ is a critical point.  Assume the formula of the theorem
is true for $r\le r'.$  Let $r''$ be the smallest critical point larger than
$r'.$  
As there are at most finitely many critical points between $r_1$ and
any $r<r_2,$ we simply need to show the formula remains valid for 
\hbox{$r'<r\le r'',$}
and the theorem follows by induction.  Indeed, as above,
\begin{align*}
	N(f,0,r) - N(f,0,r') &\;=\;\!\!\!\!\!\!\!\! 
	\sum_{\begin{array}{c}\scriptstyle |z|=r'\\
        \scriptstyle\textrm{s.t.~}f(z)=0\end{array}} 
	\!\!\!\!\!\!\!\!\log \frac{r}{r'}\\
	&\;=\;[K(f,r')-k(f,r')]\log\frac{r}{r'}\\
	&\;=\; \log|f|_r-k(f,r')\log r-\log|a_{k(f,r')}|\\
	&\;=\; \log|f|_r-\log|f|_{r'}.\qedhere
\end{align*}
\end{proof}

\smallskip
To keep the focus on the essential ideas, I will confine myself
to a discussion of one variable in these lectures.  The several variable
theory is also well-developed.  In particular, one can consider
multivariable power series or Laurent series and define $|~|_r$ similarly
to what was done in today's lecture. One of the contributions in
\cite{CherryYeTrans} is a several variable Poisson-Jensen formula which
shows that also in several variables $\log|f|_r$ measures the quantity
of zeros of $f$ in an appropriate sense. Only power series are discussed
in \cite{CherryYeTrans}, but one can also work with 
multivariable Laurent series, as discussed in \cite{CherryRu}.
Some additional discussion of several variables, particularly in postive
characteristic, is contained in \cite{CherryToropu}.

\section{Non-Archimedean Value Distribution Theory}
\mymark{Non-Archimedean Value Distribution Theory}

In this lecture we introduce the non-Archimedean analog
of Nevanlinna's theory of value distribution.

\subsection{Nevanlinna's Theory of Value Distribution}

In a deep and beautiful theory, Nevanlinna developed quantitative
analogs of the Fundamental Theorem of Algebra for meromorphic functions.
His theory continues to be
an indispensable tool in, for example, the study of complex
dynamics and the study of meromorphic solutions to differential equations;
see \textit{e.g.} \cite{BergweilerBull} and \cite{Laine}.

\smallskip
Given a meromorphic function $f,$ given a value $a$ in
$\mathbf{P}^1(\mathbf{C}),$ and given a radius $r,$ Nevanlinna
introduced the following functions. The \textbf{counting function}
$N(f,a,r)$
counts, as a logarithmic average, the number of times $f$ takes
on the value $a$ in the disc of radius $r,$ and is precisely defined by
$$
        N(f,a,r) = \!\!\!\!\!\!\sum_{\begin{array}{c}\scriptstyle0< |z|<r\\
        \scriptstyle\textrm{s.t.~}f(z)=a\end{array}} \!\!\!\!\!\!
	\log \frac{r}{|z|}\;+\;m\log r,
$$
where values $z$ in the sum over $f(z)=a$ are repeated, according
to their multiplicity, and 
where $m$ is the order of vanishing of $f(z)-a$ at the origin,
or of $1/f(z)$ if $a=\infty.$ One can also define the
\textbf{truncated counting function} $N^{(1)}(f,a,r),$ which counts
the zeros without regard to multiplicity.

The \textbf{proximity function}
$m(f,a,r)$ measures how close the function $f$ stays to $a$ on the
circle of radius $r$ and is defined by
$$
	m(f,a,r) = \int_0^{2\pi}\log^+\left|\frac{1}{f(re^{i\theta})-a}
	\right|\,\frac{d\theta}{2\pi},
$$
where $\log^+ x$ denotes $\max\{0,\log x\}.$ If $a=\infty,$
then 
$$
	m(f,\infty,r)=\int_0^{2\pi}\log^+|f(re^{i\theta})|\,
	\frac{d\theta}{2\pi}.
$$
Nevanlinna's \textbf{characteristic function} $T(f,a,r)=m(f,a,r)+N(f,a,r)$
is the sum of the counting and proximity functions.
Nevanlinna then proved two main theorems.

\begin{theorem}[Nevanlinna's First Main Theorem]\label{NevFMT}
If $f$ is a non-constant meromorphic function on $\mathbf{C}$
and $a$ is a point in $\mathbf{C},$ then $T(f,a,r)-T(f,\infty,r)$
remains bounded as $r\to\infty.$
\end{theorem}

\begin{remark*}Nevanlinna's First Main Theorem says that
\hbox{$m(f,a,r)+N(f,a,r)$}
is, up to a bounded term, independent of $a$ as $r\to\infty.$
This has two consequences. First, since $m(f,a,r)\ge0,$ it gives
an upper-bound on the frequency with which $f$ can take on the value $a.$
In this sense, the First Main Theorem is an analog of the fact that 
a polynomial of degree $d$ takes on the value $a$ at most $d$ times.
The second consequence is that if $f$ takes on a value $a$ relatively
rarely, then it must compensate for it by remaining close to the value
$a$ on a large proportion of each circle centered at the origin.
Considering the function $e^z$ is instructive.  For values $a$ other
than $0$ and $\infty,$ $e^z,$ being periodic, takes on each value
with the same frequency.  However, the values $0$ and $\infty$ are omitted
entirely.  On the other hand, $e^z$ is close to 0 and close to $\infty$
on about half of each circle centered at the origin, whereas it is
only rarely close to any non-zero value $a.$
\end{remark*}

\begin{theorem}[Nevanlinna's Second Main Theorem]\label{NevSMT}
Let $a_1,\dots,a_q$ be $q$ distinct points in $\mathbf{P}^1(\mathbf{C}),$
and let $f$ be a non-constant meromorphic function on $\mathbf{C}.$
Then,
$$
	(q-2)T(f,\infty,r) - \sum_{j=1}^qN^{(1)}(f,a_j,r) \le
	O(\log T(f,r))
$$
as $r\to\infty$ outside an exceptional set of radii of finite Lebesgue
measure.
\end{theorem}

\begin{remark*} When $q>2,$ the inequality in the Second Main Theorem
gives a lower bound
on the sum of the counting functions, so $f$ cannot take on too many
values with lower than expected frequency.
\end{remark*}

Detailed introductions to Nevanlinna's theory can be found
in \cite{Nevanlinna}, \cite{Hayman}, and \cite{CherryYeBook}.

\smallskip
This lecture will introduce the analog of Nevanlinna's theory
in non-Archimedean function theory.  Perhaps the earliest work with
the spirit of non-Archimedean Nevanlinna theory was the work
of Adams and Straus \cite{AdamsStraus}. H\`a Huy Kho\'ai, one of
the organizers of this school, was the first to set out to systematically
develop a complete analog of Nevanlinna's theory for non-Archimedean
meromorphic functions \cite{KhoaiDuke}.

\subsection{Prescribing Zeros to Analytic Functions}

The following theorem is a well-known consequence of the Mittag-Leffler
theorem and is often covered in a first course in complex analysis.

\begin{theorem}[Mittag-Leffler]\label{MittagLeffler}
Let $D$ be a domain in $\mathbf{C},$ let $z_n$ be a discrete sequence
of distinct points
in $D,$ and let $m_n\ge1$ be a sequence of positive integers.
Then, there exists an analytic function on $D$ such that for each $n,$
$f$ has a zero of multiplicity $m_n$ at $z_n$ and such that $f$ has
no other zeros.
\end{theorem}

We now explore some non-Archimedean analogs of this theorem.
First, if $f$ is analytic in $A[r_1,r_2],$ then by
Corollary~\ref{IdentityCor}, $f$ has only finitely many zeros
in $A[r_1,r_2].$ Conversely, by taking for instance a polynomial,
given a finite number of points in $A[r_1,r_2]$ and associated multiplicities,
we can construct a function $f$ analytic on $A[r_1,r_2]$ with the
prescribed zeros and multiplicities.

\begin{theorem}\label{NonArchML}
Let $z_n$ be a sequence of distinct non-zero numbers in $\mathbf{F}$ such that
$|z_n|\to\infty.$ Let $m_n$ be a sequence of positive integers,
and let $m_0$ be a non-negative integer. Then,
$$
	f(z)=z^{m_0}\prod_{n=1}^\infty \left(1-\frac{z}{z_n}\right)^{m_n}
$$
converges to an analytic function $f$ on $\mathbf{F}$ with a zero at
$0$ of multiplicity $m_0,$ with zeros at $z_n$ with multiplicities $m_n,$
and with no other zeros.
\end{theorem}

\begin{lemma}\label{infprod}
Let $f_n$ be analytic functions on $A[r,r].$ Then,
$\prod f_n$ converges if and only if $\lim f_n=1.$
\end{lemma}

\begin{proof}
Consider the partial products,
$$
	P_N = \prod_{n=1}^N f_n.
$$
We need to check that $|P_N - P_M|_r$ tends to 
zero as $\min \lbrace N,M \rbrace \to \infty.$
Without loss of generality, assume that
$N \ge M.$  Then,
$$
	|P_M-P_N|_r = \left|\prod_{n=1}^M f_n\right|_r
	\left|1 - \prod_{n=M+1}^N f_n\right|_r.
$$
Because $|1-f_n|_r \to 0,$ there exists an $N_0$ such that for
all $n \ge N_0, |f_n|_r = 1.$
Therefore for $M \ge N_0$,
$$
	\left|\prod_{n=1}^M f_n\right|_r = 
	\left|\prod_{n=1}^{N_0} f_n\right|_r.
$$
On the other hand, if $M \ge N_0,$ then
$$
	\left|1 - \prod_{n=M+1}^N f_n\right|_r
	= \left|1 - \prod_{n=M+1}^N (1 - (1-f_n))\right|_r
	\le \sup_{n > M} |1-f_n|_r.
$$
The right hand side tends to zero as $M \to \infty$ by assumption,
and the proposition follows.
\end{proof}

\begin{proof}[Proof of Theorem~\ref{NonArchML}.]
The product in the statement of the theorem converges to an
entire function by Lemma~\ref{infprod}. The product clearly
has the prescribed zeros with the prescribed multiplicities.
Fix $r>0.$ By Proposition~\ref{keqKprop}, the product over all
$z_n$ with $|z_n|> r$ converges to a unit in $\mathcal{A}[r].$
Hence, $f$ only has the zeros prescribed by the product.
\end{proof}

The situation for finite non-bordered discs is more delicate.
A complete non-Archimedean field is called \textbf{maximally complete}
or \textbf{spherically complete} if every collection of embedded
discs $D_{i+1}\subset D_i$ has non-empty intersection.

\begin{exercise} The field $\mathbf{C}_p$ is not maximally complete.
\end{exercise}

\begin{theorem}[Lazard]\label{LazardThm}
The following are equivalent:
\begin{quote}\begin{enumerate}
\item[(a)] The field $\mathbf{F}$ is maximally complete.
\item[(b)] Given $R>0,$ given a sequence of distinct points
$z_n$ in $\mathbf{B}_{<R}$ such that $|z_n|\to R,$ and given
positive integers $m_n,$ there exists an analytic function $f$ in
$\mathcal{A}(R)$ with a zero at each $z_n$ with multiplicity $m_n$
and no other zeros.
\end{enumerate}\end{quote}
\end{theorem}

The proof of Theorem~\ref{LazardThm} will not be discussed here.
See \cite{Lazard}.

\subsection{Nevanlinna Functions and the First Main Theorem}

As in the classical complex case, the non-Archimedean First Main Theorem
is simply the Poisson-Jensen formula dressed up in new notation.

\smallskip
The non-Archimedean counting functions are defined in exactly the
same way as in the complex case. If $f$ is a meromorphic function,
the proximity function is defined by
$$
	m(f,a,r) = \log^+\left|\frac{1}{f-a}\right|_r \qquad\textrm{and}
	\qquad m(f,\infty,r)=\log^+|f|_r.
$$
The characteristic function can then be defined, just as over the complex
numbers, by $$T(f,a,r)=m(f,a,r)+N(f,a,r).$$

\begin{theorem}[First Main Theorem]\label{FMT}
If $f$ is a non-constant meromorphic function on $\mathbf{F},$  then
$T(f,a,r)-T(f,\infty,r)$ remains bounded as $r\to\infty.$
\end{theorem}

\begin{proof}
We first treat the case that $a=0.$ By Theorem~\ref{NonArchML},
we can write $f=g/h,$ where $g$ and $h$ are entire without common
zeros. The Poisson-Jensen Formula (Theorem~\ref{PoissonJensen})
then tells us that
$$
	N(f,\infty,r)=N(h,0,r)=\log|h|_r +O(1) \qquad\textrm{and}\qquad
	N(f,0,r)=N(g,0,r)=\log|g|_r + O(1),
$$
and so
$$
	N(f,\infty,r)-N(f,0,r) = \log|h|_r-\log|g|_r + O(1).
$$
Now,
$$
	m(f,\infty,r)=\max\{0,\log|g|_r-\log|h|_r\} \qquad\textrm{and}\qquad
	m(f,0,r)=\max\{0,\log|h|_r-\log|g|_r\}.
$$
Thus,
$$
	m(f,\infty,r)-m(f,0,r)= \log|g|_r-\log|h|_r,
$$
which proves the theorem when $a=0.$

\smallskip
If $a\ne0,$ then $N(f,a,r)=N(f-a,0,r)$ and $N(f,\infty,r)=N(f-a,\infty,r).$
Also,
$$
	m(f,a,r)=m(f-a,0,r) \qquad\textrm{and}\qquad
	m(f,\infty,r)=m(f-a,\infty,r)+O(1).\qedhere
$$
\end{proof}

\subsection{Hasse Derivatives}
We now examine Hasse derivatives more closely. First note that if
$n$ is a non-negative integer and $k$ is an integer not in the
interval $[0,n-1],$ then the binomial coefficient
$$
	\binom{k}{n} = \frac{k(k-1)(k-2)\cdots(k-n+1)}{n!}
$$
is well-defined and non-zero.  Thus, we can extend the notion of
Hasse derivatives to Laurent series as follows. If
$$
	f(z)=\sum_{k=-\infty}^\infty a_kz^k
	\qquad\textrm{then define}\qquad
	D^n\!f(z)=\sum_{k\in\mathbf{Z}\setminus[0,n-1]}\binom{k}{n}a_kz^{k-n},
$$
and again if $\mathbf{F}$ has characteristic zero, then
$D^n\!f = f^{(n)}/n!.$ By defining $\binom{k}{n}$ to be zero if
\hbox{$0\le k \le n-1,$} we can write this simply as
$$
	D^n\!f(z)=\sum_{k=-\infty}^\infty\binom{k}{n}a_kz^{k-n}.
$$

\begin{prop}\label{Hasseprops}
The Hasse derivatives of analytic funcions on annuli
satisfy the following basic properties:
\begin{quote}
\begin{enumerate}
\item[\textup{(i)}] $\displaystyle D^k[f+g]=D^k\!f + D^k\!g;$
\item[\textup{(ii)}] $\displaystyle D^k[fg] = \sum_{i+j=k}D^i\!fD^j\!g;$
\item[\textup{(iii)}] $\displaystyle D^i\!D^j\!f = \binom{i+j}{j}D^{i+j}\!f.$
\item[\textup{(iv)}] If $\mathbf{F}$ has positive
characteristic $p$ and $s\ge0$ is an integer, then
$\ds   D^{p^s}\!f^{p^s}= (D^1\!f)^{p^s}.$
\end{enumerate}\end{quote}
\end{prop}

\begin{proof}
Property~(i) is obvious.  To check property~(ii), write out both sides
and compare like powers of $z.$ What is needed for equality is that
for $\ell$ and $m$ integers and $i$ and $j$ non-negative integers,
one has
$$
    \sum_{i+j=k}\binom{\ell}{i}\binom{m}{j}
    =\binom{\ell+m}{k},
$$
which is nothing other than Vandermonde's Identity.
To check property~(iii), one needs the elementary identity
$$
    \binom{i+j}{j}\binom{k}{i+j}
    =\binom{k}{j}\binom{k-j}{i}.
$$
What one needs for~(iv) is that fact that for any integer $j,$
$$
    \binom{jp^s}{p^s}\equiv j \textrm{~mod~}p,
$$
which follows immediately from Lucas's Theorem.
\end{proof}

Property~(ii) in Proposition~\ref{Hasseprops} allows us to inductively
extend $D^n\!f$ to meromorphic functions $f.$ For example,
$$
	D^1\!(f)=D^1\left[\frac{f}{g}g\right]
	=gD^1\!\left(\frac{f}{g}\right)+\frac{f}{g}D^1\!g,
$$
and hence
$$
	D^1\left(\frac{f}{g}\right)=\frac{gD^1\!f-fD^1\!g}{g^2}
$$
as expected.

\subsection{Logarithmic Derivative Lemma}

Nevanlinna's first proof of his Second Main Theorem is based on 
a deep property of logarithmic derivatives, namely that if $f$
is a meromorphic function on $\mathbf{C},$ then $m(f'/f,\infty,r)$
is small relative to $T(f,\infty,r)$ outside a small exceptional set
of radii $r;$ we will not give the precise
formulation here for the complex numbers.

\smallskip
The non-Archimedean analog of the Logarithmic Derivative Lemma
is a trivial, but useful, observation.

\begin{lemma}[Logarithmic Derivative Lemma]
\label{LogDerivLemma}
Let $f$ be a 
non-Archimedean meromorphic function on the annulus
$\{z : r_1 \le |z| \le r_2\}.$  Then, for $r_1\le r \le r_2,$
$$
        \left | \frac{D^nf}{f} \right |_r \le \frac{1}{r^n}.
$$
Here $D^nf$ denotes the $n$-th Hasse derivative of $f.$
\end{lemma}

\begin{proof}
We first prove the theorem in the case that $f$ is analytic.
Write $f$ as a Laurent series $f(z) = \sum c_k z^k.$
Then,
$$
        D^n\!f(z) 
	= \sum_{k\in\Z}c_k
	\binom{k}{n}z^{k-n}.
$$
Because the binomial coefficients are integers,
they have absolute value $\le 1,$
and so
\begin{equation*}
        |D^nf|_r = \sup_{k\in\Z} |c_k| \left|\binom{k}{n}\right|
         r^{k-n} \le
        \frac{\ds\sup_{k\in\Z} |c_k|r^k}{r^n} = \frac{|f|_r}{r^n},
\end{equation*}
and hence the lemma is proven for analytic $f.$

\smallskip
We prove the lemma for meromorphic $f$ by induction on $n,$ the case
that $n=0$ being trivial.  Write $f=g/h,$ where $g$ and $h$ are analytic.
Using property~(ii) of Proposition~\ref{Hasseprops} for the Hasse derivative
extended to meromorphic functions, we find that
$$
	\frac{D^{n+1}\!f}{f} 
		= \frac{D^{n+1}(g/h)}{g/h} = \frac{D^{n+1}\!g}{g}
	-\frac{D^nf}{f}\cdot\frac{D^1\!h}{h} - \cdots -
	\frac{D^1\!f}{f}\cdot\frac{D^n\!h}{h}-\frac{D^{n+1}\!h}{h},
$$
and so the lemma follows by induction, by the analytic case proven above,
and by the fact that $|~|_r$ is a non-Archimedean absolute value.
\end{proof}

\subsubsection*{An Application}

We recall here a cute argument, sometimes atributed to Gauss, that
if $f,g$ are two complex polynomials in $\mathbf{C}[z],$ not both
constant such that $f^{-1}(0)=g^{-1}(0)$ and such that $f^{-1}(1)=g^{-1}(1),$
then $f=g.$ To prove this, assume $\deg f \ge \deg g,$ and consider
the rational function
$$
	h=\frac{f'(f-g)}{f(f-1)}.
$$
Now, the degree of the denominator is twice the degree of $f$ and
the degree of the numerator is strictly less than twice the degree of $f.$
On the other hand, if $z_0$ were a pole of $h,$ then $f(z_0)=0$ or $1,$
and so by assumption $f-g$ vanishes at $z_0.$ The presence of $f'$ in
the numerator ensures that $h$ has no multiple poles, and thus $h$ has no poles
at all. Since the degree of the denominator is larger than that of the
numerator, this means $h$ is identically zero.  Since $f$ is non-constant,
$f'$ is not identically zero, and so we must have $f=g.$
Lemma~\ref{LogDerivLemma} allows us to make this same argument for
non-Archimedean entire functions.

\begin{prop}[Adams \& Straus]\label{UniequeEntireProp}
Let $f$ and $g$ be non-Archimedean
entire functions on $\mathbf{F},$ not both constant and if 
$\mathbf{F}$ has positive characteristic $p$ such that neither $f$ nor $g$
is a pure $p$-th power. Let $a$ and $b$ be two distinct points in $\mathbf{F}.$
If $f^{-1}(a)=g^{-1}(a)$ and $f^{-1}(b)=g^{-1}(b),$ then $f=g.$
\end{prop}

\begin{proof}
Without loss of generality, assume there is a sequence of 
$r_n\to\infty$ such that $|f|_{r_n}\ge|g|_{r_n}.$ The hypotheses
imply that $|f|_{r_n}\to\infty.$ As in the polynomial
case, let
$$
	h = \frac{f'(f-g)}{(f-a)(f-b)}.
$$
Then, $h$ has no poles, and is therefore entire. On the other hand,
$$
	|h|_{r_n} = \left|\frac{(f-a)'}{f-a}\right|_{r_n}
	\cdot\frac{|f-g|_{r_n}}{|f-b|_{r_n}} \le \frac{1}{r_n},
$$
for $r_n$ large enough that $|f|_{r_n}>|b|,$ and hence $h\equiv0.$
\end{proof}

\begin{remark*} Lemma~\ref{LogDerivLemma} implies
$\log|f'/f|_r \le -\log r.$ When $r\ge1,$ the right-hand-side is bounded
as $r\to\infty.$ In the early works on non-Archimedean Nevanlinna theory,
people, including me, sometimes replaced the $-\log r$ term with $O(1).$
Kho\'ai and Tu's work \cite{KhoaiTu} made clear that keeping the $-\log r$
term is essential for applications such as the above.
\end{remark*}

\subsection{Second Main Theorem}

\subsubsection*{Second Main Theorem without Ramification}

We stated Nevanlinna's Second Main Theorem over the complex numbers
as an inequality involving
$$
	\sum_{j=1}^qN^{(1)}(f,a_j,r),
$$
a sum of truncated counting functions. Of course this inequality
implies the weaker inequality where the above sum is replaced
by
$$
	\sum_{j=1}^qN(f,a_j,r),
$$
the sum of non-truncated counting functions.  This weaker inequality over
the complex numbers is still 
much deeper than Nevanlinna's
First Main Theorem.  However, in the non-Archimedean case, Ru \cite{Ru}
made the significant observation that the non-Archimedean Second Main
Theorem without ramification (or truncation) is a simple consequence
of the First Main Theorem.

\begin{theorem}[Second Main Theorem without Ramification]
\label{RuSMT}
Let $a_1,\dots,a_q$ be $q$ distinct points in $\mathbf{P}^1(\mathbf{F}).$
Then,
$$
	(q-1)T(f,\infty,r) - \sum_{j=1}^q N(f,a_j,r) \le O(1)
$$
as $r\to\infty.$
\end{theorem}

Before proceeding to the proof, several comments are in order. First observe
that in the non-Archimedean case we have $(q-1)$ on the left hand side
rather than $(q-2)$ as in the complex case.  Thus, the stronger non-Archimedean
Picard theorem follows: namely, a non-Archimedean meromorphic function on
$\mathbf{F}$ can omit at most one value in $\mathbf{P}^1(\mathbf{F}).$
Second, there is no need for an exceptional set of radii $r$ in the
non-Archimedean case.

\begin{proof}
First assume all the $a_i$ are finite and let
$d=\min\limits_{i\ne j}|a_i-a_j|>0.$ Then, given $i\ne j$ and $r>0,$
$$
	d \le |a_j-a_i| = |a_j-a_i|_r =|(f-a_i)-(f-a_j)|_r
	\le \max\{|f-a_i|_r,|f-a_j|_r\}.
$$
This implies that given $r>0,$ there is at most one index $j_0$ such that
$|f-a_{j_0}|_r < d.$ Thus,
$$
	\sum_{j\ne j_0} [T(f,a_j,r) - N(f,a_j,r)]
	= \sum_{j\ne j_0} m(f,a_j,r)
	=\sum_{j\ne j_0} \log^+|f-a_j|_r^{-1} \le (q-1)\log^+(1/d).
$$
By the First Main Theorem, up to a bounded term, we can replace
$$
	\sum_{j\ne j_0} T(f,a_j,r)
$$
with $(q-1)T(f,\infty,r).$ Since $N(f,a_{j_0},r)\ge0$ for $r\ge1,$
we can subtract $N(f,a_{j_0},r)$ from the left-hand side to get the
theorem in the case that none of the $a_j$ are infinite.
\smallskip
If one of the $a_j=\infty,$ let $c$ be a point of $\mathbf{P}^1(\mathbf{F})$
that is not among the $a_j.$ Consider the meromorphic function
$g=(f-c)^{-1}$ and let $b_j=(a_j-c)^{-1}.$ Then, the $b_j$ are all finite
and we can apply what we've already proven to $g$ and the $b_j.$
Of course $N(g,b_j,r)=N(f,a_j,r),$ and by the First Main Theorem,
$T(f,\infty,r)=T(g,\infty,r)+O(1).$
\end{proof}

\subsubsection*{The Defect Relation}
The \textbf{defect} of a point $a$ in $\mathbf{P}^1(\mathbf{F})$ with respect
to a meromorphic function $f$ is defined by
$$
	\delta_f(a)=\liminf_{r\to\infty}\frac{m(f,a,r)}{T(f,\infty,r)}
	= 1-\limsup_{r\to\infty}\frac{N(f,a,r)}{T(f,\infty,r)}.
$$
The Second Main Theorem immediately implies:

\begin{cor}\label{RuDefect}
$\ds \sum_{a\in\mathbf{P}^1(\mathbf{F})} \delta_f(a) \le 1.
$
\end{cor}

In fact, something much stronger is true.

\begin{prop}\label{defectprop}
Given a non-constant meromorphic function $f$ on $\mathbf{F},$ there is at most
one point $a$ in $\mathbf{P}^1(\mathbf{F})$ such that $\delta_f(a)>0.$
\end{prop}

\begin{proof}
By making a projective change of coordinates if necessary, without
loss of generality assume, $\delta_f(0)>0$ and $\delta_f(\infty)>0.$
Then, $\delta_f(\infty)>0$ and $\delta_f(0)>0$ imply that
$m(f,0,r)$ and $m(f,\infty,r)$ are both positive for all sufficiently 
large $r.$ But this exactly means that both $|f|_r>1$ and $|1/f|_r>1,$
which is clearly absurd.
\end{proof}

As $\delta_f(a)\le1,$ Proposition~\ref{defectprop} is stronger
than Corollary~\ref{RuDefect}. Given $a$ in $\mathbf{P}^1(\mathbf{F})$
and $0<\delta\le1,$ there exists a meromorphic function $f$ on $\mathbf{F}$
such that $\delta_f(a)=\delta;$ see, \textit{e.g.} \cite{CherryYeTrans}.

\subsubsection*{The Second Main Theorem with Ramification}

Define,
$$
	N_{\mathrm{Ram}}(f,r) = N(f',0,r)+2N(f,\infty,r)-N(f',\infty,r).
$$
Observe that $N_{\mathrm{Ram}}(f,r)$ exactly counts the ramification
points of $f$ with multiplicity, and that by Poisson-Jensen, we can also
write
\begin{equation}\label{RamPoissonJensenEqn}
	N_{\mathrm{Ram}}(f,r) = 2N(f,\infty,r)+\log|f'|_r + O(1).
\end{equation}

\begin{theorem}[Second Main Theorem with Ramification]\label{SMT}
Let $f$ be meromorphic on $\mathbf{F}$ and assume $f'\not\equiv0.$
Let $a_1,\dots,a_q$ be $q$ distinct points in $\mathbf{P}^1(\mathbf{F}).$
Then,
{\small$$
	(q-2)T(f,\infty,r)-\sum_{j=1}^q N^{(1)}(f,a_j,r) \le
	(q-2)T(f,\infty,r)-\sum_{j=1}^q N(f,a_j,r)+N_{\mathrm{Ram}}(f,r) 
	\le -\log r + O(1).
$$}
\end{theorem}

\begin{proof}
The first inequality is clear, so we show the second.
Let $q'$ denote the number of finite $a_j,$ and if $q'<q,$ without loss
of generality, assume $a_q=\infty.$
As in the proof of Theorem~\ref{RuSMT}, let
$d=\min\limits_{i\ne j}|a_i-a_j|,$ where the minimum is
taken over the finite $a_i$ and $a_j.$ For the moment, fix $r>0.$
Again, we see that there is an index $j_0,$ which may depend on $r,$
such that for all $j\ne j_0$ and $j\le q',$ we have
$|f-a_j|_r \ge d.$
Thus,
\begin{align*}
	(q'-1)m(f,\infty,r)=(q'-1)\log^+|f|_r
	&\le \sum_{j\ne j_0} \log|f-a_j|_r + (q'-1)\log^+\frac{1}{d}+
	(q'-1)\max_{1\le j \le q'} |a_j| \\
	&= \sum_{j\ne j_0} \log|f-a_j|_r +O(1),
\end{align*}
where the $O(1)$ term is independent of $r$ and $j_0,$ and the sum
over $j\ne j_0$ means the sum over all indices $\le q'$ and different from
$j_0.$ Now,
$$
	\sum_{j\ne j_0}\log|f-a_j|_r \le
	\sum_{j=1}^{q'} \log|f-a_j|_r 
	-\log|f'|_r+\log\left|\frac{f'}{f-a_{j_0}}
	\right|_r \le \sum_{j=1}^{q'}\log|f-a_j|_r - \log|f'|_r - \log r
$$
by the Logarithmic Derivative Lemma. Hence,
$$
	(q'-1)m(f,\infty,r)\le \sum_{j=1}^{q'} \log|f-a_j|_r -\log|f'|_r -
	\log r + O(1).
$$
Since the right-hand-side does not depend on $j_0,$ we no longer need
to regard $r$ as fixed. Now we apply Poisson-Jensen to get
$$
	(q'-1)m(f,\infty,r) \le \sum_{j=1}^{q'} N(f,a_j,r) - q'N(f,\infty,r)
	+2N(f,\infty,r) -N_{\mathrm{Ram}}(f,r) - \log r + O(1).
$$
Hence,
$$
	(q'-1)T(f,\infty,r) - \sum_{j=1}^{q'} N(f,a_j,r) - N(f,\infty,r)
	+N_{\mathrm{Ram}}(f,r) \le - \log r + O(1).
$$
This is precisely the statement of the theorem when $a_q=\infty.$
When $q=q',$ the theorem follows by replacing $N(f,\infty,r)$
on the left with the larger $T(f,\infty,r).$
\end{proof}

\medskip
The Second Main Theorem was first proven 
in characteristic zero by C.~Corrales-Rodrig\'a\~nez
\cite{CorralesThesis}, but her work only became easily accessible in the
literature some years later \cite{CorralesPolonici}.  The Second Main 
Theorem was also proven independently by Kho\'ai and Quang \cite{KhoaiQuang}
and by Boutabaa \cite{BoTh}. Of course the same proof works in positive
characteristic, provided $f$ is not a $p$-th power.  The formulation
above sufficies for all applications I know of. One can state an inequality
valid for all non-constant functions in positive characteristic that
essentially amounts to the following:

\begin{cor} Let $\mathbf{F}$ have positive characteristic, let
$f$ be meromorphic on $\mathbf{F}$ such that $f'\not\equiv0,$ 
let $a_1,\dots,a_q$ be $q$ distinct points in $\mathbf{P}^1(\mathbf{F}),$
and let $s$ be a non-negative integer.
Then,
$$
	(q-2)T(f^{p^s},\infty,r) - \sum_{j=1}^qN(f^{p^s},a^{p^s},r)
	+p^sN_{\mathrm{Ram}}(f,r) \le -p^s\log r + O(1).
$$
\end{cor}

\noindent
This observation, together with the fact that $N_{\mathrm{Ram}}(f,r)$
cancels any contribution to the counting functions
$N(f^{p^s},a^{p^s},r)$ coming from points whose multiplicity
is divisible by $p^{s+1},$
is essentially the content of \cite{BoutabaaEscassut}.

\subsubsection*{The ABC Inequality}

\begin{cor}[ABC]\label{ABC}
Let $f+g=h$ be relatively prime entire functions, not all of whose
derivatives vanish identically. Then,
$$
	\log\max\{|f|_r,|g|_r,|h|_r\} \le N^{(1)}(fgh,0,r)-\log r+O(1).
$$
\end{cor}

\begin{proof}
Let $F=f/h.$ By the relatively prime assumption,
$$
	N^{(1)}(fgh,0,r)=N^{(1)}(F,0,r)+N^{(1)}(F,\infty,r)+N^{(1)}(F,1,r).
$$
Applying the Second Main Theorem,
$$
	N^{(1)}(fgh,0,r)-\log r +O(1) 
	\ge T(F,\infty,r) = \log^+\left|\frac{f}{h}
	\right|_r + N(h,0,r).
$$
Now, $N(h,0,r)=\log|h|_r+O(1)$ by Poisson-Jensen, so
$$
	\log^+\left|\frac{f}{h}\right|_r +N(h,0,r)
	= \log\max\{|f|_r,|h|_r\} + O(1).
$$
Since $g=h-f,$ we know $|g|_r \le \max\{|f|_r,|h|_r\},$
and the corollary follows.
\end{proof}

Hu and Yang undertook a systematic study of generalized ABC-inequalities for
non-Archimedean entire functions: \cite{HuYangABCone}, \cite{HuYangABCtwo},
and \cite{HuYangABCthree}. See \cite{CherryToropu} for positive 
characteristic.

\subsubsection*{The Defect Relation Again}

The \textbf{ramification defect} $\theta_f(a)$
of a point $a$ in $\mathbf{P}^1(\mathbf{F})$
with respect to a meromorphic function $f$ is defined by
$$
	\theta_f(a)=1-\limsup_{r\to\infty}\frac{N^{(1)}(f,a,r)}{T(f,\infty,r)}.
$$
The Second Main Theorem immediately implies

\begin{cor} If $f$ is a meromorphic function on $\mathbf{F}$ such
that $f'\not\equiv0,$ then 
$$
	\sum_{a \in \mathbf{P}^1(\mathbf{F})} \theta_f(a) \le 2,
$$
with strict inequality if $f$ is a rational function.
\end{cor}

Unlike the case with ordinary defects, it is possible for 
$\theta_f(a)$ to be positive for more than one value of $a.$
Not much is known about ramification defects for non-Archimedean
meromorphic functions, and it would be interesting to say anything
non-trivial about them.

\subsubsection*{Some Applications}

A value $a$ is called \textbf{totally ramified} for a meromorphic
function of $f$ if every point in $f^{-1}(a)$ is a ramification point.
For example, $1$ and $-1$ are totally ramified values for the sine and
cosine functions.

\begin{cor}If $f$ is a non-Archimedean meromorphic function such
that $f'\not\equiv0,$ then $f$ has at most three totally ramified 
values.
\end{cor}

\begin{proof}
Suppose $a_1,\dots,a_4$ are totally ramified values.
Then, the Second Main Theorem says:
$$
	2T(f,\infty,r) - \sum_{j=1}^4 N^{(1)}(f,a_j,r) \le -\log r + O(1).
$$
But because the $a_j$ are totally ramified,
$$
	N^{(1)}(f,a_j,r)\le \frac{1}{2}N(f,a_j,r) \le
\frac{1}{2}T(f,\infty,r)+O(1),
$$
where the second inequality follows from the First Main Theorem.
We thus conclude \hbox{$\log r \le O(1)$} and thereby reach a contradiction.
\end{proof}

It is easy to check that $0,$ $1,$ and $\infty$ are totally ramified
values for the rational function
$$
	f(z)=\frac{(z^2-1)^2}{(z^2+1)^2},
$$
and thus the theorem cannot be improved.  Over the complex numbers,
the Weierstrass $\wp$ function has four totally ramified values;
that this is the most possible is a consequence of the Second Main
Theorem.

\begin{exercise}\label{EntireTotRamif}
Show that a non-Archimedean entire function $f$
such that $f'\not\equiv0$ can have at most one finite totally ramified
value. What happens if the hypothesis $f'\not\equiv0$ is dropped
in positive characteristic?
\end{exercise}

\begin{theorem}[Adams \& Straus]\label{meroshare}
Let $f$ and $g$ be two meromorphic functions on $\mathbf{F}.$ If
$\mathbf{F}$ has characteristic zero, assume $f$ and $g$ are not both
constant. If $\mathbf{F}$ has positive characteristic, assume that
neither $f'\equiv0$ nor $g'\equiv0.$ Let $a_1,\dots,a_4$ be distinct
elements of $\mathbf{P}^1(\mathbf{F})$ and assume $f^{-1}(a_j)=g^{-1}(a_j)$
for $j=1,\dots,4.$ Then, $f=g.$
\end{theorem}

\begin{proof}
We treat the case that $f'\not\equiv0$ and $g'\not\equiv0.$
It isn't difficult to modify the proof to allow one of the functions
to be constant.
Without loss of generality, assume none of the $a_j$ are infinity.
It is easy to see that $T(f-g,\infty,r)\le T(f,\infty,r)+T(g,\infty,r).$
By hypothesis,
$$
	N^{(1)}(f-g,0,r) \ge \sum_{j=1}^4N^{(1)}(f,a_j,r) = \sum_{j=1}^4
	N^{(1)}(g,a_j,r).
$$
Applying the Second Main Theorem to both $f$ and $g$ we conclude that
\begin{align*}
	2T(f,\infty,r)+2T(g,\infty,r) &\le 2\sum_{j=1}^4N^{(1)}(f,a_j,r)
	-2\log r + O(1)\\
	&\le 2N^{(1)}(f-g,0,r)-2\log r+O(1) \\
	&\le 2T(f-g,0,r)-2\log r + O(1)\\
	&\le 2T(f,\infty,r)+2T(g,\infty,r)-2\log r + O(1),
\end{align*}
which is a contradiction.
\end{proof}

The example
$$
	f(z)=\frac{z}{z^2-z+1} \qquad\textrm{and}\qquad
	g(z)=\frac{z^2}{z^2-z+1}
$$
shows that Theorem~\ref{meroshare} is best possible since
$$
	f^{-1}(0)=g^{-1}(0), \qquad f^{-1}(1)=g^{-1}(1),
	\qquad\textrm{and} \qquad f^{-1}(\infty)=g^{-1}(\infty).
$$

\subsection{An's Defect Relation}

Ru's observation that non-Archimedean inequalities of Second Main Theorem
type follow from the First Main Theorem has been quite important.
As an example, Ta Thi Ho\`ai An \cite{An} proved the following:

\begin{theorem}[An's Defect Relation]\label{AnThm}
Let $X\subset\mathbf{P}^N$ be a projective variety and let
\hbox{$f:\mathbf{F}\rightarrow X$}
be a non-constant non-Archimedean analytic map.
Let $D_1,\dots,D_q$ be hypersurfaces in $\mathbf{P}^N$ in general
position with $X,$ and assume that the image of $f$ is not
completely contained in any of the $D_j.$ Then,
$$
	\sum_{j=1}^q \delta_f(D_j) \le \dim X.
$$
\end{theorem}

\begin{remark*} An proved Theorem~\ref{AnThm} while she was visiting
ICTP as a Junior Associate.
\end{remark*}

\noindent
Here, in general position with $X$ means that the intersection of
any $\dim X+1$ of the $D_i$ and $X$ is empty. Here the defects
$\delta_f(D_j)$ measure $f$ encountering the hypersurface $D_i$
with less than expected frequency. In particular, 
a non-constant $f$ cannot omit more than $\dim X$ of the $D_i,$
unless it is completely contained in one of them. A defect relation
such as this is unique to non-Archimedean analysis and has no
counterpart in complex value distribution theory, in the sense that the
dimension bounds the defect sum and that the bound is derived from the
First Main Theorem.  The deeper  defect inequality of
Eremenko and Sodin \cite{EremenkoSodin} over the complex
numbers takes a similar form to
An's inequality, but with $\dim X$ replaced
by $2\dim X,$ and it is true for entirely different reasons.

\subsection{Concluding Remarks}

The Second Main Theorem, with ramification, can also be proven for
maps encountering hyperplanes in projective space; see \cite{BoCrv},
\cite{KhoaiTu}, and \cite{CherryYeTrans}. The techniques of this
section can be used to prove that any non-Archimedean analytic map from
the affine line $\mathbf{A}^1$ to an algebraic curve of positive genus
must be constant; see \cite{CherryWang} for details. This was first
proven by Berkovich \cite{Berkovich} using his theory of analytic
spaces. My lecture during the workshop in the third week will be
about the degeneracy of images of non-Archimedean analytic maps
to projective varities omitting divisors with sufficiently many components.
Perhaps one of the most interesting things to investigate in 
non-Archimedean function theory is analogs of Big Picard theorems.
For instance, one can prove \cite{CherryBigPic} that a non-Archimedean
analytic map from a punctured disc to an elliptic curve with good
reduction must always extend across the puncture.  But this need not
be true for elliptic curves with bad reduction.   That whether a Big 
Picard type theorem is true or not can depend on the reduction type
of the target is a phenomenon completely foreign to the complex analytic
situation.

\section{Benedetto's Island Theorems}
\mymark{Benedetto's Island Theorems}

\subsection{Ahlfors Theory of Covering Surfaces}
In work that won him one of the first Fields Medals, Ahlfors \cite{AhlforsCov}
developed a theory of covering surfaces that both gave a geometric
interpretation of Nevanlinna's Second Main Theorem as a generalization
of the Gauss-Bonet Formula and extended it to the distribution of
``domains,'' rather than ``values,'' and also to classes of mappings
more general than meromorphic functions, for instance quasiconformal
mappings. The Ahlfors Five Islands Theorem was a consequence of his
covering theory and has been an important tool in the study of complex
dynamics; see \cite{Bergweiler}.
Recall that as a consequence of the Second Main Theorem, a
meromorphic function can have at most four totally ramified values.
That means if $f$ is a meromorphic function on 
$\mathbf{C}$ and $a_1,\dots,a_5$ are five distinct values 
in $\mathbf{P}^1(\mathbf{C}),$
then there must be a point $z_0$ in $\mathbf{C}$ such that $f(z_0)$
is one of the $a_j$ and $f(z_0)=a_j$ with multiplicity one. 
The Five Islands Theorem is the same statement, but with the values
$a_j$ replaced by domains.

\begin{theorem}[Ahlfors Five Island Theorem]
Let $f$ be a meromorphic function on $\mathbf{C}$ and let
$D_1,\dots,D_5$ be five simply connected domains in $\mathbf{P}^1(\mathbf{C})$
with disjoint closures.  Then, there exists an open set $U$ in $\mathbf{C}$
such that $f$ is a conformal bijection between $U$ and one of the
domains $D_1,\dots,D_5.$
\end{theorem}

\begin{remark*} The theorem gets its name because  Ahlfors thought
of the five domains $D_1,\dots,D_5$ as islands on the Riemann sphere.
The theorem says that given five islands, a meromorphic function must
cover at least one of the islands injectively.
\end{remark*}

In two significant papers, \cite{BenedettoEntire} and \cite{BenedettoMero},
R.~L.~Benedetto investigated non-Archimedean
analogs of the Ahlfors island theorems. This lecture is intended as
an introduction to Benedetto's work.

\subsection{A Non-Archimedean Riemann Mapping Theorem?}

The islands in Ahlfors's theorems are simply connected domains.
What should the non-Archimedean analog be?

\begin{prop}\label{discimage}
Let $f$ be a non-constant analytic function on $\mathbf{B}_{\le r}$ and let
$R=|f-f(0)|_r.$ Then,
$$
	f(\mathbf{B}_{\le r}) = \{ w \in \mathbf{F} : |w-f(0)| \le R\}.
$$
\end{prop}

Proposition~\ref{discimage} says that the image of a disc under a 
non-constant non-Archimedean analytic function must be another disc.
Thus, the non-Archimedean analog of the Riemann mapping theorem would
be the trivial
statement that given any bordered disc $D$ in $\mathbf{F},$ there
is an analytic map from $D$ to $\mathbf{B}_{\le 1}.$

\begin{proof} Let $z$ in $\mathbf{B}_{\le r}.$ Then,
$$
	|f(z)-f(0)| \le |f-f(0)|_r
$$
by the Maximum Modulus Principle, and hence
$$
	f(\mathbf{B}_{\le r}) \subset \{w \in \mathbf{F} : |w-f(0)|\le R\}.
$$
Now let $w$ be such that $|w-f(0)|\le R.$ Then,
$$
	|f-w|_r = |f-f(0)+f(0)-w|_r \le \max\{|f-f(0)|_r, |w-f(0)|\}=R.
$$
Write 
$$
	f(z)=f(0)+\sum_{k=1}^\infty c_k z^k, \qquad\textrm{and so}\qquad
	f(z)-w=f(0)-w+\sum_{k=1}^\infty c_k z^k.
$$
By assumption,
$$
	\sup_{k\ge1}|c_k|r^k = R, \qquad\textrm{and so}\qquad
	K(f-w,r)\ge 1.
$$
Hence, $f(z)-w$ has a zero $\mathbf{B}_{\le r}$ by Theorem~\ref{NewtonPolyThm}.
\end{proof}

\begin{exercise} If $f$ is a non-constant analytic function on 
$\mathbf{B}_{<r},$ then $f(\mathbf{B}_{<r})$ is an unbordered disc,
including the possibility that $f(\mathbf{B}_{<r})=\mathbf{F},$
which can be viewed as an unbordered disc of infinite radius.
\end{exercise}

By a bordered disc in $\mathbf{P}^1(\mathbf{F}),$ we mean either a
bordered disc in $\mathbf{F}$ or a set of the form
$$
	\{w \in \mathbf{F} : |w|\ge R>0\}\cup\{\infty\}.
$$
By an unbordered disc in $\mathbf{P}^1(\mathbf{F}),$ we mean either
an unbordered disc in $\mathbf{F}$ (including the possiblility of 
$\mathbf{F}$ itself), or a set of the form
$$
	\{w \in \mathbf{F} : |w|>R\ge 0\}\cup\{\infty\}.
$$

\begin{exercise} If $f$ is a non-constant meromorphic function
on $\mathbf{B}_{\le r}$ (resp.\ $\mathbf{B}_{<r}$), then
$f(\mathbf{B}_{\le r})$ (resp.\ $f(\mathbf{B}_{<r})$) is
either all of $\mathbf{P}^1(\mathbf{F})$ or
a bordered (resp.\ unbordered) disc in $\mathbf{P}^1(\mathbf{F}).$
\end{exercise}

\subsection{Non-Archimedean Analogs of the theorems of Bloch, Landau, 
Schottky, and Koebe}

In \cite{BenedettoEntire}, Benedetto formulated and proved non-Archimedean
analogs of the classical theorems of Bloch, Landau, Schottky, and Koebe.
These are stated here as exercises for the reader. The first is the
most difficult and depends on the characteristics of $\mathbf{F}$ and
$\Ftilde.$ See \cite{BenedettoEntire} for solutions to these exercises,
as well as further commentary and references for the classical complex
analogs.

\begin{exercise}[Non-Archimedean Bloch's Constant]\label{Bloch}
\begin{quote}
\begin{enumerate}
\item[~]~
\item[(i)] If $\mathrm{char}\;\mathbf{F}=0,$ let
$$
	B=\left\{\begin{array}{ll}
	1 &\textrm{~if~}\mathrm{char}\;\Ftilde=0\\
	|p|^{1/(p-1)}&\textrm{~if~}\mathrm{char}\;\Ftilde=p>0.
	\end{array}\right.
$$
Let $f$ be analytic on $\mathbf{B}_{\le 1}$ normalized so that
$f(0)=0$ and $f'(0)=1.$
Then, there is an unbordered disc $U$ in $\mathbf{B}_{\le 1}$ such that
$f$ is injective on $U$ and such that $f(U)$ is an unbordered disc of
radius $B.$ Moreover, there exists an analytic function $f$ on
$\mathbf{B}_{\le 1}$ such that $f(0)=0,$ $f'(0)=1,$ and such that
if $U$ is any bordered or unbordered disc in $\mathbf{B}_{\le 1}$
on which $f$ is injective, then $f(U)$ does not contain a bordered
disc of radius $B.$
\item[(ii)] If $\mathrm{char}\;\mathbf{F}>0,$ then given $\varepsilon>0,$
there exists an analytic function $f$ on $\mathbf{B}_{\le 1}$ such
that $f(0)=0,$ such that $f'(0)=1,$ and such that if $U$ is any bordered
or unbordered disc in $\mathbf{B}_{\le 1}$ such that $f$ 
is injective on $U,$ then $f(U)$ does not
contain a bordered disc of radius $\varepsilon.$
\end{enumerate}
\end{quote}
\end{exercise}

\begin{remark*} The constant $B$ in Exercise~\ref{Bloch} (and zero
if $\mathrm{char}\;\mathbf{F}>0$) can be called
the non-Archimedean Bloch constant, and the examples of the exercise
show the constant is sharp. Over the complex numbers, the existence
of a positive constant $B$ such that if $f$ is holomorphic on the unit
disc in $\mathbf{C}$ normalized so that $f'(0)=1,$ then $f$ injectively
covers some disc of radius $B$ is a theorem of Bloch. The sharp value
of $B$ in the complex case is a long-standing conjecture, and remains
unproven.
\end{remark*}

\begin{exercise}[Non-Archimedean Landau's Constant]\label{LandauConst}
If $f$ is analytic on $\mathbf{B}_{\le 1}$ normalized such that 
$f(0)=0$ and $f'(0)=1,$ then 
\hbox{$f(\mathbf{B}_{\le 1})\supseteq\mathbf{B}_{\le 1}.$}
\end{exercise}

\begin{remark*} Exercise~\ref{LandauConst} says that the non-Archimedean
Landau constant is 1, and the value $1$ is clearly best
possible. As with Bloch's constant, the sharp value of Landau's
constant over the complex numbers is conjectured, but not yet proven.
\end{remark*}

\begin{exercise}[Non-Archimedean Koebe 1/4-Theorem]
If $f$ is analytic and injective on $\mathbf{B}_{\le 1}$ and normalized
so that $f(0)=0$ and $f'(0)=1,$ then
$f(\mathbf{B}_{\le 1})=\mathbf{B}_{\le 1}.$
\end{exercise}

\begin{exercise}[Non-Archimedean Landau Theorem]
Let $f$ be analytic and zero free on $\mathbf{B}_{\le 1}.$
Then $|f'(0)|<|f(0)|.$
\end{exercise}

\begin{exercise}[Non-Archimedean Schottky Theorem]
Let $f$ be analytic and zero free on $\mathbf{B}_{\le 1}.$
Then, $|f(z)|=|f(0)|$ for all $z$ in $\mathbf{B}_{\le 1}.$
\end{exercise}

\subsection{Island Theorems}

\subsubsection*{Entire Functions}

Given Exercise~\ref{EntireTotRamif}, one might expect the following
statement to be a non-Archimedean analog of the Ahlfors Island Theorem
for entire functions.

\begin{statement}\label{EntireIsland}
Let $D_1$ and $D_2$ be two disjoint unbordered discs in $\mathbf{F}$
and let $f$ be an entire function on $\mathbf{F}$ such that $f'\not\equiv0.$
Then, there exists an unbordered disc $U$ in $\mathbf{F}$ such that $f$
is injective on $U$ and such that $f(U)=D_1$ or $f(U)=D_2.$
\end{statement}

Unfortunately, as we will see in a moment, 
Statement~\ref{EntireIsland} is false if
$\mathrm{char}\;\Ftilde>0,$ even in the case that
$\mathrm{char}\;\mathbf{F}=0.$ We will also see that 
Statement~\ref{EntireIsland} is true if $\mathrm{char}\;\Ftilde=0.$

\medskip
Before looking at some examples in positive characteristic, we
give a general proposition.

\begin{prop}\label{InjectiveImage}
Let $f$ be analytic and injective on a bordered disc $D$ of radius $r$
containing the point $a.$
Then, $f(D)$ is a bordered disc of radius at most $r|f'(a)|.$
\end{prop}

\begin{proof}
That $f(D)$ is a bordered disc is Proposition~\ref{discimage}.
Without loss of generality, assume $f$ is given by a power series of
the form
$$
	f(z)=\sum_{n=1}^\infty a_nz^n.
$$
Let $b>r|a_1|.$ If $b$ is in $f(\mathbf{B}_{\le r}),$ then
$$
	|a_1|r < |b|\le |f|_r = \sup_{n\ge1}|a_n|r^n,
$$
and so $K(f,r)\ge2,$ and $f$ is not injective on $\mathbf{B}_{\le r}.$
\end{proof}

\subsubsection*{Difficulties in positive characteristic}

We now explain a fundamental difference between the case of analytic
functions over the complex numbers and non-Archimedean analytic functions
when $\mathrm{char}\;\Ftilde>0.$

\begin{exercise}\label{ComplexInverse}
Let $G$ be a domain in $\mathbf{C}$ and let $f$ be analytic on $G.$
Let $R$ be the set of ramification points of $f,$ \textit{i.e.,}
$$
	R=\{z \in G : f'(z)=0\}.
$$
Let $B=f(R)$ be the set of branch points in $f(G).$ Let $D$ be a
simply connected domain in $f(G)\setminus B.$ Then, there exists
an analytic function $g,$ called a branch of $f^{-1},$
on $D$ with values in $G$ such that
$f(g(z))\equiv z$ on $D.$ This implies that the open set $U=g(D)$
has the property that $f$ is injective on $U$ and $f(U)=D.$
\end{exercise}

When $\mathrm{char}\;\Ftilde=0,$ the same property holds in
the non-Archimedean case.

\begin{exercise}\label{NonArchBranch}
Let $f$ be analytic on a disc $\mathbf{B}_{\le r}.$
Let $R\subset\mathbf{B}_{\le r}$ be the set of ramification points
and $B=f(R)$ be the set of branch points. Let $D$ be a bordered or
unbordered disc in $f(\mathbf{B}_{\le r})\setminus B.$ Then, there
exists an analytic function $g$ on $D$ with values in 
$\mathbf{B}_{\le r}$ such that $f(g(z))\equiv z$ on $D.$
\end{exercise}

When $\mathrm{char}\;\mathbf{F}=p>0,$ the statement
in Exercise~\ref{NonArchBranch} is spectacularly false.
The polynomial $f(z)=z+z^p$ is such that $f'(z)\equiv1,$ so that
$f$ has no ramification or branch points, but $f$ fails to have an 
inverse function. The same phenomenon persists even when
$\mathrm{char}\;\mathbf{F}=0<p=\mathrm{char}\;\Ftilde.$
In this case, $R$ is not empty, but rather
$$
	R = \{\zeta \in \mathbf{F} : \zeta^{p-1}=-1/p\},
$$
and so $|\zeta|=|p|^{-1/(p-1)}>1$ for each $\zeta$ in $R.$
Hence, $|f(\zeta)|=|p|^{-p/(p-1)}>1,$ and so 
$\mathbf{B}_{\le 1}\cap B=\emptyset.$ Nonetheless, $f$ is not injective
on any disc mapping onto $\mathbf{B}_{\le 1}.$ [Exercise: prove this.]

\subsubsection*{An Example When $\mathrm{char}\;\mathbf{F}>0.$}

Consider the case that $\mathrm{char}\;\mathbf{F}=\mathrm{char}\;\Ftilde=p>0.$
Clearly the hypothesis $f'\not\equiv0$ in Statement~\ref{EntireIsland}
is necessary because for pure
$p$-th powers, such as $f(z)=z^p,$ every value is totally ramified and
$f$ is nowhere injective. However, Statement~\ref{EntireIsland} remains
false even with this hypothesis.

\begin{example}[\textup{\cite[pp.~598]{BenedettoEntire}}]\label{zplusczp}
Let $\mathrm{char}\;\mathbf{F}=p>0,$ let $\varepsilon>0,$ and let
$c$ be an element of $\mathbf{F}$ such that $|c|>\varepsilon^{-(p-1)}.$
Then, $f(z)=z+cz^p$ is not injective on any bordered disc of radius
$\varepsilon$ (characteristic $p$!) and hence $f$ does not injectively cover
any bordered disc of radius $\varepsilon$ by Proposition~\ref{InjectiveImage},
since $f'(z)\equiv1.$
\end{example}

Example~\ref{zplusczp} not only shows that Statement~\ref{EntireIsland}
is false, but it also shows that no island theorem can be true for all
entire functions (or even polynomials) when $\mathrm{char}\;\mathbf{F}>0,$
even if the number of islands is increased or one requires the islands to
be very small.

\subsubsection*{An Example When 
$\mathrm{char}\;\mathbf{F}=0<p=\mathrm{char}\;\Ftilde.$}

We now show that Statement~\ref{EntireIsland} is false
when $\mathrm{char}\;\Ftilde = p > 0 = \mathrm{char}\;\mathbf{F},$
even if one increases the number of islands or adds an additional
restriction that the islands be ``small.''

\begin{example} Let $\mathrm{char}\;\Ftilde=p>0=\mathrm{char}\;\mathbf{F}.$
Let $a_i$ be infinitely many points in $\mathbf{B}_{\le 1}$ such that
$|a_i|=1$ for all $i$ and such that 
$|a_i-a_j|=1$ for all $i\ne j,$ which is possible since $\Ftilde$
is algebraically closed, and hence infinite. Let $1\ge\varepsilon>0$ and
choose $n$ such that $|p|^n < \varepsilon.$ Let $D_i$ be unbordered discs
containing $a_i$ of radius $\varepsilon.$ Let $f(z)=z^{p^n}.$ 
Then, $f$ does not injectively cover any of the $D_i,$ which 
are disjoint since $|a_i-a_j|=1.$ Indeed, for each $i,$
let $\xi_i$ be a point in $\mathbf{F}$ such that
$f(\xi_i)=a_i,$ and note that $|\xi_i|=1.$ Let $U$ be a bordered disc
containing $\xi_i$ on which $f$ is injective. Clearly, $f$ is not injective
on $\mathbf{B}_{\le 1},$ and so the radius of $U$ is at most 1.
It then follows from Proposition~\ref{InjectiveImage} that $f(U)$ is a 
disc containing $a_i$ of radius at most $|p|^n < \varepsilon,$
and therefore $D_i$ is not injectively covered by $f.$
\end{example}

\subsubsection*{Benedetto's Island Theorem for Analytic Functions on a Disc}

As we have seen, Statement~\ref{EntireIsland} is false when
$\mathrm{char}\;\Ftilde>0.$ I introduced this lecture with the statements
of Ahlfors's island theorems for functions meromorphic or holomorphic
on $\mathbf{C}.$ In fact, Ahlfors's theorems apply to functions meromorphic
or holomorphic on a disc that satisfy certain additional
hypotheses. The Ahlfors
island theorems for functions meromorphic or analytic on $\mathbf{C}$
then follow by showing that if $f$ is a non-constant meromorphic or analytic
function on $\mathbf{C},$ then when $f$ is restricted to
sufficiently large discs, it satisfies the additional hypotheses of the 
associated disc island theorem. Benedetto's
point of view is that although the island theorem for non-Archimedean
entire functions
is not true in general, there is a good non-Archimedean analog of 
Ahlfors's island theorem for functions analytic on discs.

Before stating the theorem, we introduce some convenient notation.
For a non-Archimedean analytic function, define
$$
	f^\#(z) = \left \{ \begin{array}{ll}
	\ds\frac{\ds|f'(z)|}{\ds \max\{1,|f(z)|\}}
	&\textrm{~if~}f(z)\ne\infty\\
\noalign{\vskip 3pt}
	\ds \left|\left(\frac{1}{f}\right)'(z)\right|
	&\textrm{~if~}f(z)=\infty.
	\end{array}\right.
$$
The quantity $f^\#$ is a non-Archimedean analog of the spherical
derivative in complex analysis, and it is also convenient to adopt
the notation
$$
	||f^\#||_r = \frac{|f'|_r}{\max\{1,|f|_r\}}.
$$
One sees that $f^\#$ is the natural measure of the distortion of $f$
considered as a map to $\mathbf{P}^1(\mathbf{F}).$

\begin{theorem}[Benedetto's Analytic Island Theorem]\label{AnalyticIsland}
Let $D_1$ and $D_2$ be two disjoint
unbordered discs in $\mathbf{F},$ each of finite
radius. There exist explicit
constants $C_1$ and $C_2$ depending only on $D_1,$ $D_2$ and
the characteristics
of $\mathbf{F}$ and $\Ftilde$ with $C_2=0$ when $\mathrm{char}\;\Ftilde=0,$
such that the following holds. Given $f$ analytic on $\mathbf{B}_{<1}$
with $f^\#(0)>C_1$ and $r||f^\#||_r\ge C_2$ for some $0<r<1,$ then there
exists an unbordered disc $U$ in $\mathbf{B}_{<1}$ such that $f$ is
injective on $U$ and such that $f(U)=D_1$ or $f(U)=D_2.$
\end{theorem}

\begin{remark*}The hypothesis $f^\#(0)> C_1$ is a natural necessary
hypothesis for an analytic function on a disc. This hypothesis ensures
that both $D_1$ and $D_2$ are in the image of $f,$ and clearly without
some hypothesis to ensure that the islands are in the image of $f,$
no such theorem would be possible.  The hypothesis
$r||f^\#||_r\ge C_2$ for some $r$ is automatically
satisfied when $\mathrm{char}\;\Ftilde=0$ and
has the effect of ruling out the
positive characteristic pathologies that we explored above.
Benedetto's point of view is that this second hypothesis
is in some sense a natural non-Archimedean analog to
the hypothesis in Ahlfors's island theorem for meromorphic functions
on a disc that the mean covering number be sufficiently big with respect
to the relative boundary length, two notions from Ahlfors's theory
of covering surfaces that I will not attempt to make precise here.
A significant difference, though, between the complex and non-Archimedean
cases is that in the complex case, if $f$ is a non-constant
meromorphic function on $\mathbf{C},$ then $f(rz)$ will satisfy
Ahlfors's hypothesis for all $r$ sufficiently large, and therefore 
result in his island theorem for non-constant meromorphic functions
on $\mathbf{C},$ whereas in the non-Arcimedean setting when
$\mathrm{char}\;\Ftilde>0,$ we have seen examples of functions $f$
such that none of the functions $f_a=f(az),$ no matter how large $|a|$ is,
satisfy the hypothesis $r||f_a^\#||_r\ge C_2$ for some $0<r<1.$
\end{remark*}

I will not give the idea of the proof of Theorem~\ref{AnalyticIsland}
in the most interesting case when $\mathrm{char}\;\Ftilde>0$ in
these lectures and instead simply refer the reader to Benedetto's
paper \cite{BenedettoEntire}. I will instead prove a special case
of the theorem in the case that $\mathrm{char}\;\Ftilde=0$ that
explains the essential steps in this most simple of cases.

\begin{prop}Assume $\mathrm{char}\;\Ftilde=0.$ Let
$$
	D_0 = \{z \in\mathbf{F} : |z| < 1\} \qquad\textit{and}\qquad
	D_1 = \{ z \in\mathbf{F} : |z-1| < 1\}.
$$
Let $f$ be analytic on $\mathbf{B}_{\le1}$ such that 
$f^\#(0)\ge 1.$ Then, there exists an unbordered disc $U$ in 
$\mathbf{B}_{<1}$ such that $f$ is injective on $U$ and such that
$f(U)=D_0$ or $f(U)=D_1.$
\end{prop}

\begin{proof}
Write $f$ as a power series,
$$
	f(z)=\sum_{k=0}^\infty a_k z^k.
$$
The hypothesis $|f'(0)|\ge\max\{1,|f(0)|\}$ exactly says that
$|a_1|\ge|a_0-b|$ for all $|b|\le1.$ In particular,
$K(f-b,1)\ge1$ for all $|b|\le1,$ and hence
$$
	f(\mathbf{B}_{\le1})\supseteq\mathbf{B}_{\le1}.
$$
Note that I have just provided the solution to Exercise~\ref{LandauConst},
and so far the characteristic zero hypothesis has not been used.

\smallskip
I now explain what the hypothesis $\mathrm{char}\;\Ftilde=0$ provides.
Let $K=K(f,1)=K(f-1,1)\ge1.$ By definition, for $j>K,$ we have
$|a_j|<|a_K|.$ The characteristic zero hypothesis then implies
$$
	|j||a_j| < |K||a_K| \quad\textrm{for~}j>K\qquad\textrm{and}\qquad
	|K||a_K| \ge |j||a_j| \quad\textrm{for~}1\le j\le K,
$$
since $|j|=|K|=1.$ Hence, 
\begin{equation}\label{charzero}
K(f',0,1)=K-1.
\end{equation}

\medskip
For a subset $X\subseteq \mathbf{B}_{\le1},$ we define the following
counting functions. We let $n(f,\{0,1\},X)$ denote the number of 
points $z$ in $X$ such that $f(z)$ is in $\{0,1\},$ and we count with points
repeated according to multiplicity. Similarly, we let
$n(f',0,X)$ denote the number of zeros of $f'$ in $X,$ again repeated
according to multiplicity. The key fact we need to show in order to
prove the proposition is that there exists a point $z_0$ in 
$\mathbf{B}_{\le 1}$ with $f(z_0)\in\{0,1\}$ and 
such that if $U$ is any disc containing $z_0,$
then
\begin{equation}\label{islandeqn}
	n(f,\{0,1\},U)>2n(f',0,U).
\end{equation}
Since $K(f,1)=K(f-1,1)=K,$ by Theorem~\ref{WeierstrassThm},
there are $2K$ points $z,$ counting multiplicity,
in $\mathbf{B}_{\le1}$ such that $f(z)\in\{0,1\}.$
Suppose that~(\ref{islandeqn}) is false. Then, we may choose
disjoint discs $X_i$ containing each of the finitely many 
points $z_j$ such that $f(z_j)\in\{0,1\}.$ Note that we are not
claiming that each $X_i$ contains exactly one $z_j;$ the same disc
$X_i$ may contain several of the $z_j.$ Nonetheless, by the disjointness
of the $X_i,$ we would have
$$
	2K=n(f,\{0,1\},\mathbf{B}_{\le1}) =
	\sum_i n(f,\{0,1\},X_i) \le
	2\sum_i n(f',0,X_i) \le 2n(f',0,\mathbf{B}_{\le1})=2(K-1),
$$
by~(\ref{charzero}),
which would be a contradiction. Thus,~(\ref{islandeqn}) is established.
for some point $z_0.$ The inequality in~(\ref{islandeqn}) implies
that $f'(z_0)\ne0.$ Without loss of generality, assume $f(z_0)=0,$
for otherwise we may replace $f$ with $1-f.$ Also, by a linear change
of coordinate, we may assume without loss of generality that $z_0=0.$
Now let $R<1$ be the smallest radius such that $|f|_R=1.$
I claim that if $U=\mathbf{B}_{<R},$ then $f$ is injective on $U$
and $f(U)=D_0.$ That $f(U)=D_0$ follows immediately from
$|f|_R=1,$ since if $|b|<1,$ we can find $r<R$ with $|f|_r=|b|.$
Let $r<R.$ Because $|f|_r<1,$ we have that
$$
	n(f,\{0,1\},\mathbf{B}_{\le r}) = K(f,r).
$$
Also,
$$
	n(f',0,\mathbf{B}_{\le r}) = K(f',r) = K(f,r)-1,
$$
where again the last inequality follows from the hypothesis that
$\mathrm{char}\;\Ftilde=0$ in the same way as (\ref{charzero}).
By~(\ref{islandeqn}),
$$
	2(K(f,r)-1) = 2n(f',0,\mathbf{B}_{\le r}) <
		n(f,\{0,1\},\mathbf{B}_{\le r}) = K(f,r),
$$
and hence $K(f,r)<2.$ Thus, $f$ is injective on $\mathbf{B}_{\le r},$
and hence also on $U.$
\end{proof}

\subsubsection*{Benedetto's Island Theorem for Meromorphic Functions
on a Disc}

I will conclude this lecture with the statement of Benedetto's Four Island
Theorem for meromorphic functions on a disc. When $\mathrm{char}\;\Ftilde>0,$
there is another added complication for meromorphic functions beyond what
was already present for analytic functions.  Thus, it is more complicated
to formulate a hypothesis to exclude the full range of pathologies one
can meet in positive characteristic. Benedetto states his hypothesis
using the language of Berkovich analytic spaces. We will let 
$\mathcal{P}^1$ be the Berkovich analytic space associated
to $\mathbf{P}^1$ and we will let $\mathcal{B}_{<1}$ be the
Berkovich space associated to $\mathbf{B}_{<1}.$ Then a meromorphic
function $f$ on $\mathbf{B}_{<1}$ naturally extends to a mapping
from $\mathcal{B}_{<1}$ to $\mathcal{P}^1.$ Also, recall that each
$\nu$ in $\mathcal{B}_{<1}$ is naturally associated to a
multiplicative semi-norm on the ring of analytic functions on
$\mathbf{B}_{<1},$ and so therefore extends to a semi-norm 
$||~||_\nu$ on the
field of meromorphic functions on $\mathbf{B}_{<1}.$
Recall also that to each point $\nu$ in $\mathbf{B}_{<1},$ one can
associate an embedded family of bordered discs
$$
	B(a_i,r_i)=\{z\in\mathbf{F}:|z-a_i|\le r_i\} \subset
	B(a_{i-1},r_{i-1})=\{z\in\mathbf{F}: |z-a_{i-1}|\le r_{i-1}\}
	\subset \mathbf{B}_{<1},
$$
and hence one naturally associates to $\nu$ a ``radius,'' $r(\nu)$ defined
by $r(\nu)=\inf r_i.$ Although the family of discs $B(a_i,r_i)$ is not
uniquely determined by $\nu,$ the quantity $r(\nu)$ is. Finally,
for $\nu$ with $r(\nu)>0$ in $\mathcal{B}_{<1}$ and $f$ meromorphic on
$\mathbf{B}_{<1},$ we let $||f^\#||_\nu$ denote
$$
	||f^\#||_\nu = \frac{||f'||_\nu}{\max\{1,||f||_\nu\}}.
$$
We are then finally able to state Benedetto's theorem.

\begin{theorem}[Benedetto's Meromorphic Island Theorem]\label{MeromorphicIsland}
Let $D_1,\dots, D_4$ be four disjoint
unbordered discs in $\mathbf{P}^1.$
Let $\mu$ be a point of 
$\mathcal{P}^1$
such that no connected component (in the Berkovich topology)
of $\mathcal{P}^1\setminus\{\mu\}$ intersects more than
two of the $D_i.$ There exist explicit
constants $C_1$ and $C_2$ depending only on the $D_i,$ on $\mu,$ and on
the characteristics
of $\mathbf{F}$ and $\Ftilde$ with $C_2=0$ when $\mathrm{char}\;\Ftilde=0,$
such that the following holds. Given $f$ meromorphic on $\mathbf{B}_{<1}$
with $f^\#(0)>C_1$ and $r(\nu)||f^\#||_\nu\ge C_2$ for all points
$\nu$ in $\mathcal{B}_{<1}$ such that $f(\nu)=\mu,$ then there
exists an unbordered disc $U$ in $\mathbf{B}_{<1}$ such that $f$ is
injective on $U$ and such that $f(U)=D_i$ for some $i$ from $1$ to $4.$
\end{theorem}

I simply refer the reader to \cite{BenedettoMero} for a
discussion of the proof, but we conclude by mentioning one example.

\begin{example}[\textup{\cite[Ex.~6.2]{BenedettoMero}}]
Examples
of the form
$$
	\left(z+\frac{c}{z^{p^n}}\right)\prod_{i=1}^N
	\left(1+\frac{c}{(z-a_i)^{p^n}}\right),
$$
where $|a_i-a_j|=1$ if $i\ne j$ and $|a_i|=1$ for all $i$ show that
even if one allows very many very small islands when $\mathrm{char}\;\Ftilde=p,$
one cannot significantly weaken the hypotheses of
Theorem~\ref{MeromorphicIsland} as one can show that
such examples satisfy the inequality involving
$r(\nu)||f^\#||_\nu$ for some $\nu$ with $f(\nu)=\mu,$ but not for all
$\nu$ with $f(\nu)=\mu,$ and  for these functions,
the conclusion of the theorem is
false for small islands around the $a_i.$
\end{example}

\end{document}